\documentclass[%
 aip,
 amsmath,amssymb,
 reprint,%
]{revtex4-1}

\usepackage{graphicx}
\usepackage{epstopdf}
\usepackage{dcolumn}
\usepackage{bm}

\usepackage[utf8]{inputenc}
\usepackage[T1]{fontenc}
\usepackage{mathptmx}
\usepackage{float}
\usepackage{amsmath,amssymb,mathrsfs}
\numberwithin{equation}{section}

\begin{document}

\preprint{AIP/123-QED}

\title[Discovering mean residence time and escape probability]
{Discovering mean residence time and escape probability from data of stochastic dynamical systems}

\author{Dengfeng Wu}
\email{wudengfeng@mail.cncnc.edu.cn}
 \affiliation{
College of Computer Science and Technology, Changchun Normal University, Changchun 130032, China 
}%
\author{Miaomiao Fu}
 \email{mmfucaathy@163.com}
\affiliation{
College of Mathematics, Changchun Normal University, Changchun 130032,  China 
}%

\author{Jinqiao Duan}
\email{duan@iit.edu}
\affiliation{%
Department of Applied Mathematics, Illinois Institute of Technology, Chicago 60616, U.S.A. 
}%

\date{\today}

\begin{abstract}
We present a method to learn mean residence time and escape probability from data modeled
 by stochastic    differential equations.  This method is a combination of machine
learning from data (to extract stochastic differential equations as models) and stochastic dynamics (to quantify dynamical behaviors with deterministic tools).
The goal is to learn and  understand stochastic dynamics based on data.
This method is applicable to sample path data  collected  from  complex systems, as long as these systems    can  be modeled as stochastic differential equations.
\end{abstract}

\maketitle

\begin{quotation}
Stochastic dynamical systems are   appropriate models for randomly influenced systems.
Understanding the complex dynamical behaviors of these systems is a challenge in   diverse areas of science and engineering.
In deterministic dynamical systems, invariant manifolds and
 other invariant structures provide global information for dynamical evolution.
For stochastic dynamical systems,  better quantitative analysis and understanding is needed because of
limitations in our current analytical skills or computation capability.
Fortunately, researchers are increasingly using data-driven methods for system identification and the discovery of dynamics.
In this work, we propose a new approach to determine some   computable dynamical quantities from data, such as the
mean residence time and escape probability, which offer insights into global dynamics under uncertainty.
We demonstrate the algorithm to be effective and robust, by reproducing known dynamics and evaluating errors for   several  prototypical stochastic dynamical systems with Brownian motions.

\end{quotation}

\section{\label{sec:level1}Introduction}
\setcounter{equation}{0}
\renewcommand\theequation{1.\arabic{equation}}

Stochastic dynamical systems  arise in modeling   molecular dynamics,
 mechanical and electrical engineering, climate dynamics, geophysical and environmental systems, among others.
Advances in machine learning and data science are leading to new progresses  in the analysis and
understanding of complex   dynamics for systems with
  massive observation data sets. Despite the rapid
development of tools to extract governing equations from  data, there has been slow
progress in distilling  quantities that may be used to explore stochastic dynamics.
It is desirable  to extract  deterministic quantities
that carry dynamical information of stochastic differential equations (SDEs).
These deterministic quantities include moments for solution paths, probability density
functions for solution paths, mean residence time and escape probability \cite{Duan02, Duan03}.
These concepts help us to understand various phenomena in complex systems under uncertainty \cite{Arnold98}.

In this present paper, we will present a new approach to extract the  underlying
deterministic quantities,   mean residence time and escape probability,  that describe  certain aspects of  stochastic dynamics.
In fact, mean residence time for a stochastic  dynamical system quantifies how long the system
stays in a region, and escape probability describes the likelihood of  a system transition
    from one regime to another. Fortunately, these deterministic quantities
can be determined  by solving an elliptic partial differential equation as   in Duan \cite{Duan03},
once the underlying stochastic differential equation model for the system evolution is discovered from data.

We assume that a data set is composed of sample paths governed by  a  stochastic differential equation (SDE)   in $\mathbb{R}^n$,
\begin{eqnarray}\label{dX}
dX_{t}=b(X_{t})dt+\sigma(X_{t})dB_{t},
\end{eqnarray}
where $b$ is an n-dimensional vector function, $\sigma$ is an $n\times m$ matrix function, and
$B_{t}$ is an m-dimensional Brownian motion.
This is the customary, probabilistic way of writing the equation
$$
\frac{dX_{t}}{dt}=b(X_{t})  +\sigma(X_{t}) \frac{dB_{t}}{dt}.
$$
Often, $b$ is called `drift' and $\sigma$ is called `diffusion'.  Assume that  $b$ and $\sigma$
      satisfy an appropriate local Lipschitz condition  as follows:
\begin{eqnarray*}
\|b(x)-b(y)\|+ \|\sigma(x)-\sigma(y)\|\leq K_{N}\|x-y\|,
\end{eqnarray*}
for $\|x\|\leq N$, $\|y\|\leq N$ and  $N>0$. Here the Lipschitz constant $K_{N}$ depends on the positive number $N
$. The generator for this SDE system is a linear second-order differential operator:
\begin{eqnarray}\label{A}
Ag=b\cdot(\bigtriangledown g)+\frac{1}{2}\mathrm{Tr}[\sigma \sigma^{T}H(g)],~~~~g\in\mathbf{H}^{2}_{0}(\mathbb{R}^n),
\end{eqnarray}
where $H$ denotes the Hessian matrix of a multivariate function and $\mathrm{Tr}$  denotes the trace of a matrix.
We view $A$ as a linear differential operator in Hilbert space $\mathbf{L}^{2}(\mathbb{R}^n)$ with domain of
definition $D(A)=\mathbf{H}^{2}_{0}(\mathbb{R}^n)$.  Furthermore, we assume
that the generator $A$ is uniformly elliptic.
That is for $x\in D$ and all $\xi \in \mathbb{R}^n$,
there exists a positive constant $C$ such that
\begin{eqnarray}\label{A1}
\sum^{n}_{i,j=1}(\sigma(x)\sigma^{T}(x))_{i,j}\xi_{i}\xi_{j}\geq C|\xi|^2.
\end{eqnarray}

We discuss mean residence time and escape probability as deterministic quantities
that carry dynamical information for solution orbits of \eqref{dX}.
For a bounded domain $D \subset \mathbb{R}^{n}$ (with boundary $\partial D$),
the first exit time for a solution orbit starting at $x\in D$ is a stopping time defined as
\begin{eqnarray*}
\tau_{D}(\omega)\triangleq
\inf \{t>0:X_{0}=x,X_{t}\in \partial D\}.
\end{eqnarray*}

The mean residence time is defined as $u(x)\triangleq \mathbb{E}\tau_{D}(\omega)$.
It is the mean residence time of a particle initially at $x$
inside $D$ until the particle first hits the boundary $\partial D$ or escapes from $D$.
 It turns out  that the mean residence time $u$ of stochastic system \eqref{dX}
can be determined by solving a deterministic partial differential equation as follows:
\begin{eqnarray}
Au &=& -1,\label{u1}\\
u|_{\partial D} &=& 0,\label{u2}
\end{eqnarray}
where $A$ is the generator defined as \eqref{A}. To prove the existence and uniqueness of $u$,
we recall the Dynkin's formula (Theorem 7.4.1 in Oksendal\cite{Oksendal}):
\begin{eqnarray*}
\mathbf{E} ^x [f(X_{\tau_D})]=f(x)+\mathbf{E} ^x[\int_{0}^{\tau_D}Af(X_s)ds],
\end{eqnarray*}
for $f$ in the domain of definition of the generator $A$, $\mathbf{E} ^x$ is the expectation with respect to
the probability law $Q^x$   induced by a solution process $X_t$ starting at $x\in \mathbb{R}^n$.
We take the continuous, bounded boundary value $\phi =0$ and the continuous inhomogeneous term $g=1$.
By Theorem 9.3.3 in Oksendal\cite{Oksendal}, which is a consequence of Dynkin's formula, we know that the linear
expectation $\mathbf{E} ^x[\phi(X_{\tau_{D}})]+\mathbf{E} ^x[\int_{0}^{\tau_D}g(X_s)ds]$, which is just
$\mathbf{E}_{\tau_D} = u(x)$, solves $Au=-1$ with the boundary condition $u|_{\partial D}=0$.
Thus we can numerically  compute mean residence time $u$ by solving elliptic partial differential equations,
so we know how long the stochastic system \eqref{dX} stays in the region $D$. For more details, see  Duan\cite{Duan03}.

The escape probability is the likelihood that an orbit starting inside a domain $D$,
exits from this domain first through a specific part $\Gamma$ of the boundary $\partial D$.
Let $\Gamma$ be a subset of the boundary $\partial D$. We define the escape probability $p(x)$ from $D$ through $\Gamma$ as
the likelihood that $X_t$ starting at $x$ exits from $D$ first through
$p(x)=\mathbb{P}\{X_{\tau_{\partial D}}\in \Gamma\}$. We will show that the escape probability $p(x)$ solves a linear
elliptic partial differential equation, with a specifically chosen  Dirichlet condition as follow:
\begin{eqnarray}
Ap &=& 0,\label{p1}\\
p|_{\Gamma} &=& 1,\label{p2}\\
p|_{\partial D \setminus \Gamma} &=& 0,\label{p3}
\end{eqnarray}
where $A$ is the generator defined as \eqref{A}. Taking
\begin{equation*}
\phi (x)=\left\{
\begin{aligned}
1, & & {x\in \Gamma,}\\
0, & & {x\in \partial D\setminus \Gamma,}
\end{aligned}
\right.
\end{equation*}
we have
\begin{eqnarray*}
\mathbf{E}[\phi(X_{\tau_{\partial D}}(x))]
&=&\int_{\{\omega:X_{\tau_{\partial D}}(x)\in \Gamma\}}\phi (X_{\tau_{\partial D}}(x))d\mathbb{P}(\omega)\\
&&+\int_{\{\omega:X_{\tau_{\partial D}}(x)\in \partial D\setminus \Gamma\}}\phi (X_{\tau_{\partial D}}(x))d\mathbb{P}(\omega)\\
&=& \mathbb{P}\{\omega: X_{\tau_{\partial D}}(x)\in \Gamma \} \\
&=& p(x).
\end{eqnarray*}
This means that $\mathbf{E}[\phi(X_{\tau_{\partial D}}(x))]$ is the escape probability $p(x)$ that we are looking for.
We know $Ap=0$ together with $p|_{\Gamma}=1$ and $p|_{\partial D\setminus \Gamma}=0$ by Theorem 9.2.14 in Okendal\cite{Oksendal}.  For more details, see  Duan\cite{Duan03}.
Moreover we suppose that solution orbits, i.e., `particles'  are initially uniformly distributed in $D$,
then the average escape probability $P$ that a trajectory will leave $D$ through $\Gamma$ is $P=\frac{1}{D}\int_{D}p(x)dx$.

To obtain the mean residence time $u$ and escape probability $p$ from sample path data  $X_{True}$   (with the underlying model  \eqref{dX}), we propose the following machine learning  algorithm:
First we collect sample path data $X_{True}$  by sample-wisely simulating \eqref{dX} via Euler method (and treat these data as our observation data), then  we try to learn the   stochastic dynamical system  model  from these data  with the following model ansatz
\begin{eqnarray} \label{learnmodel}
\dot{X}=\Theta \cdot  \Xi,
\end{eqnarray}
where the basis $\Theta$ consists of polynomial   functions $\{1,x,y,\cdots,x^{N}y^{N}\}$  together with the (generalized) time derivative of Brownian motion $dB_{t}/dt$, and $\Xi$ is the coefficient (or coordinate, or weight) under this basis. Here $\cdot$ denotes scalar product. Each column $\xi^k=[\xi_1^k, \xi_2^k, \cdots, \xi_N^k]$ of $\Xi$ is a vector of coefficients, determining
which terms are active in the right-hand side for one of the row  equations in \eqref{dX}.
Brunton et al. \cite{Kutz} considered such a learned model,  with a basis not containing noise (i.e., $dB_{t}/dt$).

Furthermore, we set up a regression problem
 to determine coefficient $\Xi=[\xi^1,\cdots,\xi^n]$ by minimizing the mean-square discrepancy (other metrics are possible)  between the  data  $X_{True}$ (many samples),  and the solution $X_{Learn}$ (many samples)  for the learned (i.e., extracted)  governing model \eqref{learnmodel}.
This will provide us the \emph{learned} drift $b$ and \emph{learned} diffusion $\sigma$, and thus we also have the extracted stochastic    model \eqref{learnmodel}, together with   the \emph{learned} generator $A$.

Finally,  with this learned generator $A$,  we   compute mean residence time  $u$ and escape probability  $p$,        by solving the  deterministic partial differential equations \eqref{u1}-\eqref{u2} and  \eqref{p1}-\eqref{p3}, respectively.
We verify that our algorithm is effective by estimating the error  of  the maximal mean residence time\\
$$\mbox{error}(u)=\max_{ D}\|u_{Learn}-u_{True}\|,   $$
  and the error of  the average escape probability \\
$$\mbox{error}(P)=\|P_{Learn}-P_{True}\|,   $$
 between the learned system   \eqref{learnmodel} and  the original system  \eqref{dX}.

This  paper is organized as follows.
In section 2,  we apply our method to learn    a two-dimensional quasigeostrophic meandering jet model with
  additive noise and multiplicative noise,  and compute  the mean residence time and escape probability.
In section 3, we illustrate our method to learn    two three-dimensional systems  (a  linear damped oscillator and the well-known  Lorenz system).  We summarize and conclude in section 4 .

\section{\label{sec:level1} Learning two dimensional stochastic dynamical systems}

\setcounter{equation}{0}
\renewcommand\theequation{2.\arabic{equation}}

Zonal shear flows occur naturally in both oceans and the atmosphere. The
Gulf stream is a well known example. Based on RAFOS floats observations of the Gulf Steam,
Bower  et  al. \cite{Bower01, Bower02} viewed the fluid motion as a steady and eastward meander
propagation in the moving frame,  and divided the velocity field into three  regimes: a central jet,
exterior retrograde motion, and intermediate closed circulations above meander troughs and below crests.
There was no exchange occurring between the three regimes in this model. Fluid particles in the
intermediate regime execute periodic motion but never escape.
However, float observations  reveal  that  several particle  trajectories pass
through meander crests and troughs and then leave jet, which indicate that exchange
does occur across some part of the Gulf Stream.  Samelson \cite{Samelson01} modified the
basic model and considered three different types of variability: a time-dependent
spatially uniform meridional velocity superimposed on the basic flow, a time-dependent meander amplitude,
and a propagating plane wave superimposed on the basic flow. In particular, he only considered time-periodic
variabilities. Subsequently, Samelson\cite{Samelson02}, del-Castillo-Negrete et al.\cite{Del}, Pratt et al.\cite{Pratt01},
Beigie et al.\cite{Wiggins}, Duan et al.\cite{Duan01} and Brannan et al.\cite{Duan02}  have obtained a series of outcomes from the point of view of dynamical systems.

\subsection{\label{sec:level2}Discovering a kinematical model for a  two-dimensional meandering jet}

\begin{figure*}[htb]
\includegraphics[scale=0.20]{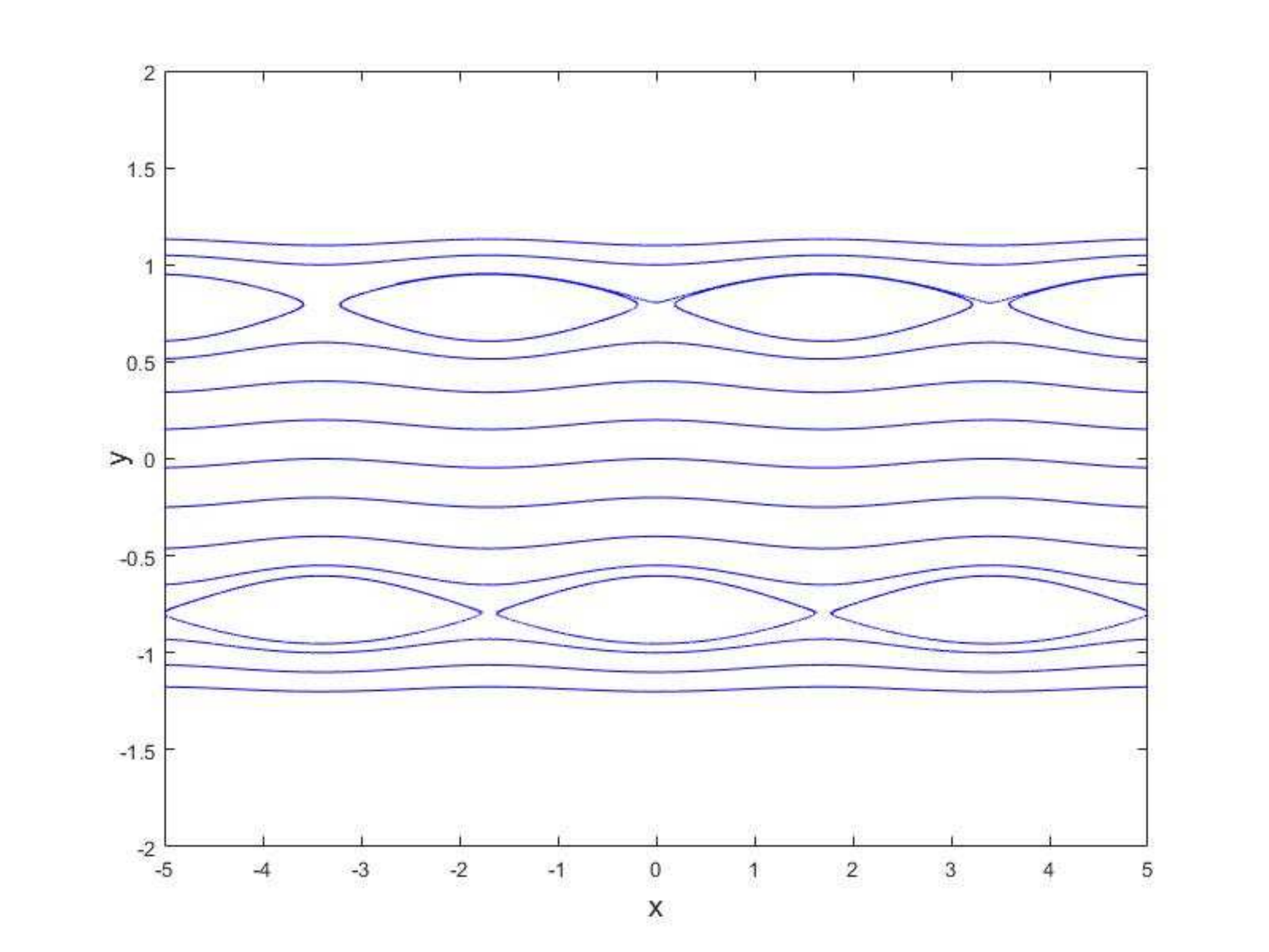}
\includegraphics[scale=0.20]{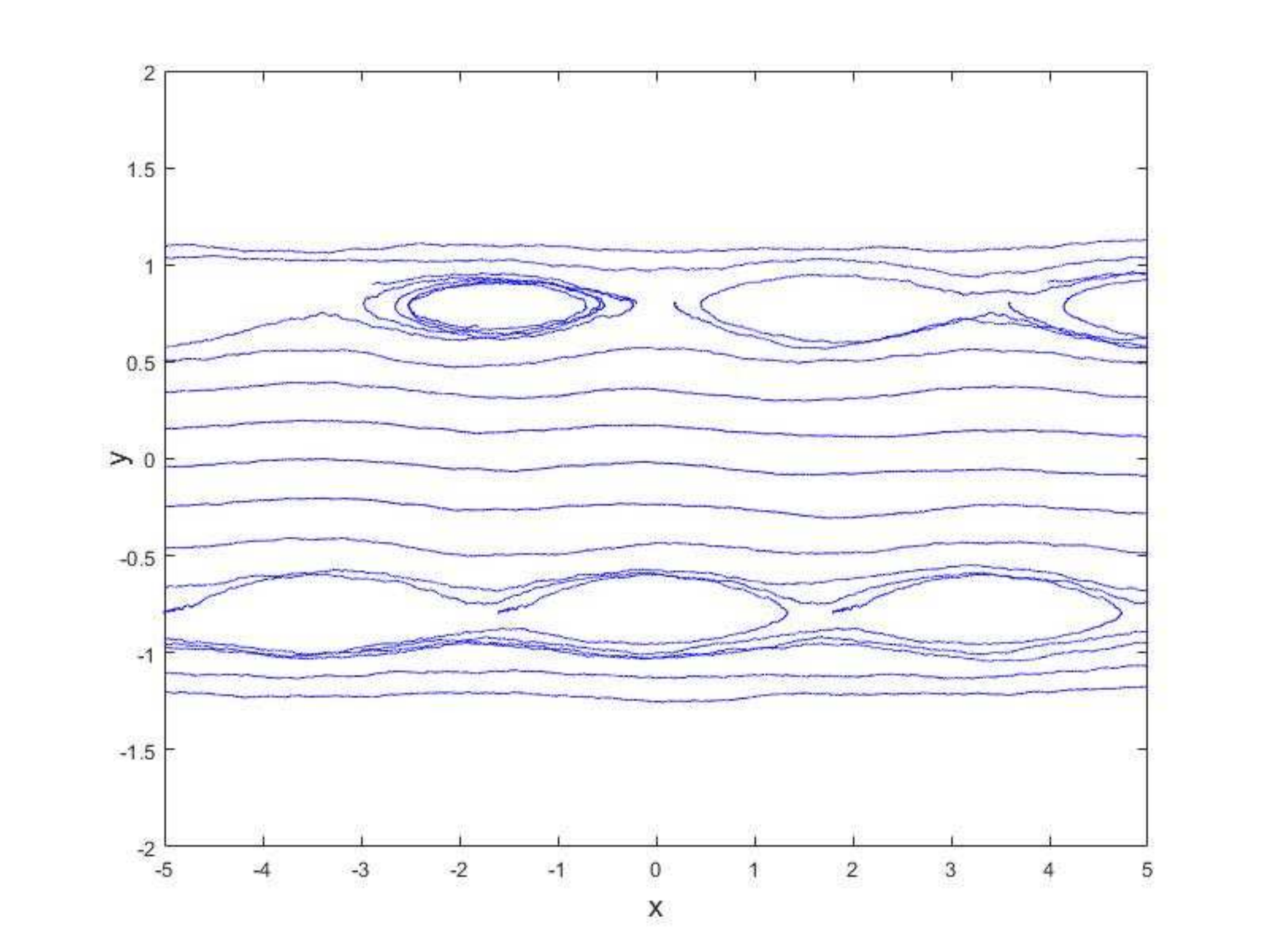}
\includegraphics[scale=0.20]{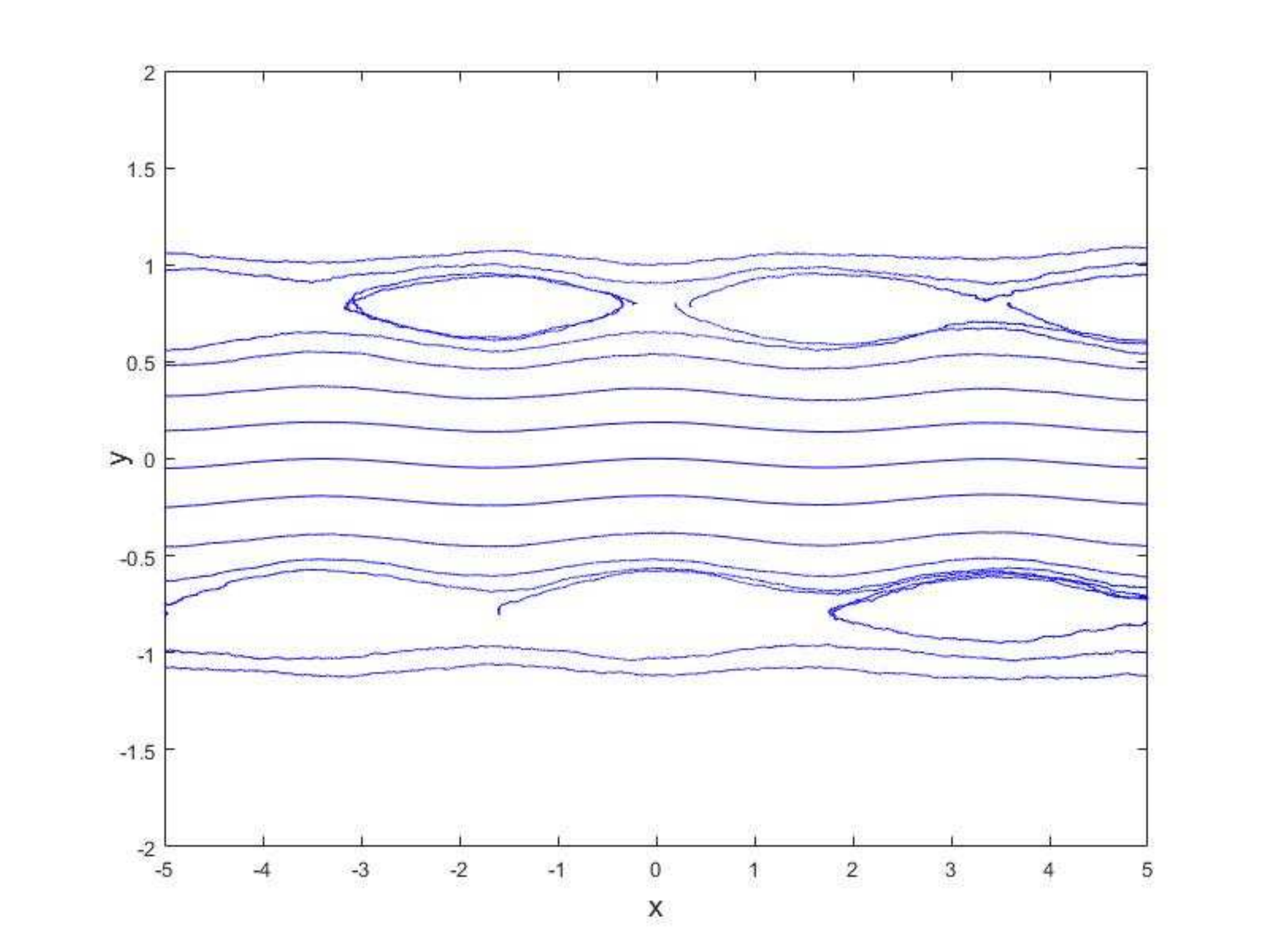}
\caption{Lagrangian fluid trajectories of the original  meandering jet  system  \eqref{ad1}-\eqref{mu2} with 20 initial positions: The fluid motion in the deterministic basic  model ($\sigma=0$) (left)
is steady, with no exchange of fluid between the eddies and the jet stream.
Under additive noise ($\sigma=\sqrt{0.3}$) (middle) and multiplicative noise ($\sigma=\sqrt{0.3}$) (right),
there are   fluid exchanges between   eddies and the jet stream.}
\label{figure 1}
\end{figure*}

As an approximate solution to the quasigeostrophic   model for two-dimensional
geophysical flows \cite{Ped}, the basic Bickley  jet \cite{Del} is
\begin{equation*}
\psi (x,y)=-\tanh (y)+a sech ^2(y)\cos (kx) +cy,
\end{equation*}
where $a=0.01$, $c=\frac{1}{3}(1+\sqrt{1-\frac{3}{2}\beta})$, $k=\sqrt{6c}$, $0\leq\beta\leq \frac{2}{3}$.
In this paper, we take $\beta=\frac13$.
Note $\beta=\frac{2\Omega}{r} \cos \theta$ is the meridional derivative of the Coriolis parameter,
where $\Omega$ is the rotation rate of the earth and $r$ is the earth's radius and $\theta$ is the latitude.
This stream function $\psi(x,y)$ defines  the basic meandering jet   system
 $\dot{x} = -\psi_{y}, \;\; \dot{y} = \psi_{x}$.
 The phase portrait for this deterministic system is in
 FIG. \ref{figure 1} (left). As it shows no fluid exchange between the eddies and the jet, it is not an appropriate model for a meandering jet.

Thus we incorporate  random wind forcing and other fluctuations, either additive or multiplicative, in the model for the meandering jet and
consequently, we consider the following two stochastic dynamical  systems
\begin{eqnarray}
dx &=& -\psi_{y}dt+\sigma dB_1,\label{ad1}\\
dy &=& \psi_{x}dt+\sigma dB_2,\label{ad2}
\end{eqnarray}
and
\begin{eqnarray}
dx &=& -\psi_{y}dt+\sigma x dB_1,\label{mu1}\\
dy &=& \psi_{x}dt+\sigma y dB_2,\label{mu2}
\end{eqnarray}
where the noise intensity $\sigma$ satisfies $0<\sigma<1$, and $B_1(t)$,  $B_2(t)$ are two independent Brownian motions.
For both stochastic systems, a number of sample solution orbits are shown in  FIG. \ref{figure 1} (middle, right).
The central  jet   and two rows of recirculation eddies,
which are called the northern and southern recirculation regions, are still visible.
Outside the recirculation regimes are the exterior retrograde regimes.
There is no exchange between regimes in unperturbed model,
but exchange does occur in other two models with additive or multiplicative noises. These two stochastic  models are more appropriate for modeling a meandering jet (as a simplified model for
  the Gulf Stream), as observations indicate fluid exchange.


A wide range of dynamical systems with intrinsic noise may be modeled
with stochastic differential equations.
However, identifying an appropriate  SDE model from intermittent observations
of the system is challenging, particularly if the dynamical process is
nonlinear and the observations are noisy and indirect \cite{GaoDuan}.
Brunton et al. \cite{Kutz},  Zhang et al. \cite{Lin},   and Dunker et al. \cite{Bohner} presented
some approaches to discover or learn   governing physical laws  from data,  with underlying   differential
equations models.



We use  our machine learning  algorithm  to   data,  from  system \eqref{ad1}-\eqref{ad2}
or system  \eqref{mu1}- \eqref{mu2},
to learn their governing laws in terms of  the following model
\begin{eqnarray}
\dot{x} &=& \Theta \xi^1,  \label{learn1} \\
\dot{y} &=& \Theta \xi^2,    \label{learn2}
\end{eqnarray}
where $\Theta$ is a set of basis functions and $\Xi=[\xi^1,\xi^2]$ defines the  coefficients (or weights). To achieve this, the idea is as follows:
First we collect a time series of the system state $(x(t_i), y(t_i)), i=1,\cdots, 20000,$ from system \eqref{ad1}-\eqref{ad2}
or system  \eqref{mu1}- \eqref{mu2} with noise intensity $\sigma=\sqrt{0.3}$,   numerically by the Euler  method (as our observation data
$(x^{i}_{True}, y^{i}_{True})$).
Then we construct a library $\Theta$  consisting of basis functions  with
polynomials $\{1,x,y,\cdots,y^5\}$ and time derivative of Brownian motions  $\{dB_1/dt,dB_2/dt\}$.
  We determine each column of coefficients $\xi^k=[\xi_1^k,\xi_2^k,\dots,\xi_{23}^k],k=1,2,$   by  minimizing discrepancies
$\mathbf{E}\sum_{i=1}^{20000}\|x^{i}_{True}-x^{i}_{Learn}\|^2$ and
$\mathbf{E}\sum_{i=1}^{20000}\|y^{i}_{True}-y^{i}_{Learn}\|^2$.  Here $(x^{i}_{Learn}, y^{i}_{Learn})$
is the numerical solution of the learned meandering jet model  \eqref{learn1}- \eqref{learn2} by Euler method.
 See Table \ref{tab:11}  and  Table \ref{tab:12}.

\begin{table}[ht]
\caption{Identified coefficients for  the  learned meandering jet system  \eqref{learn1}-\eqref{learn2}.
Data is numerically generated via Euler method for  \eqref{ad1}-\eqref{ad2} (additive noise) with initial position $(-0.2,0.8)^{T}$ in the eddy,  for time $t\in [0,200]$ with stepsize $0.01$.}
\label{tab:11}
\begin{ruledtabular}
\begin{tabular}{ccc}
basis & $\dot{x}$ & $\dot{y}$\\
\noalign{\smallskip}\hline\noalign{\smallskip}
$1$ &  2.2114650e+00  & 1.1356876e+00 \\
$x$ &  9.1670929e-01 & -3.1053406e-01 \\
$y$ &   -8.5337156e+00 & -7.8830089e+00 \\
$x^2$ &    2.1581933e-01 & -5.8882267e-02 \\
$xy$ &    -3.5406793e+00 &  1.2086885e+00\\
$y^2$ &   1.4806990e+01 &  2.1509332e+01 \\
$x^3$ &   2.3163030e-02 & -2.0687376e-03 \\
$x^2 y$ &    -6.4404598e-01 &  1.6556400e-02 \\
$xy^2$ &   5.2163263e+00 & -2.3046831e+00 \\
$y^3$ &  -1.3568792e+01 & -2.9159973e+01 \\
$x^4$ & 6.6869672e-04 &  1.1181427e-02 \\
$x^3 y$ &  -7.3808685e-02  & 6.4410700e-02 \\
$x^2 y^2$ &  5.2977754e-01 &  2.2568117e-01 \\
$xy^3$ &   -3.5783954e+00 &  2.2945214e+00\\
$y^4$ & 5.6314933e+00 &  1.9776092e+01 \\
$x^5$ &   2.8510400e-04 &  1.2382566e-03 \\
$x^4 y$ & -8.1958722e-04 & -8.0288516e-04\\
$x^3 y^2$ &   4.2275033e-02 & -3.6886874e-02\\
$x^2 y^3$ &    -1.1615968e-01 & -1.7364994e-01 \\
$xy^4$ &   9.9449264e-01 & -9.0959557e-01\\
$y^5$ &  -6.8100788e-01 & -5.3876952e+00\\
$dB1/dt$ &  5.4773030e-01 & -9.7392788e-06 \\
$dB2/dt$ &   2.7076551e-07  & 5.4771837e-01 \\
\end{tabular}
\end{ruledtabular}
\end{table}

\begin{table}[ht]
\caption{Identified coefficients for  the  learned meandering jet system  \eqref{learn1}-\eqref{learn2}.
Data is numerically generated via Euler method for  \eqref{mu1}-\eqref{mu2} (multiplicative  noise) with initial position $(-0.2,0.8)^{T}$ in the eddy,   for time $t\in [0,200]$ with stepsize $0.01$.}
\label{tab:12}
\begin{ruledtabular}
\begin{tabular}{ccc}
basis & $\dot{x}$ & $\dot{y}$\\
\noalign{\smallskip}\hline\noalign{\smallskip}
$1$     & 4.6791584e+00 & -7.5867935e+00\\
$x$     &6.3006784e-01  &-2.1087326e+00 \\
$y$     & -2.5483586e+01  & 4.3335141e+01\\
$x^2$   & 2.6109263e-01 & -1.5952411e-01 \\
$xy$    &  -1.7420741e+00 &  1.0564395e+01\\
$y^2$   & 6.1169058e+01 & -9.4722221e+01\\
$x^3$   &  -3.1824190e-03 & -1.0056540e-01\\
$x^2 y$ &  -9.6615354e-01 & -2.4673718e-01\\
$xy^2$  & 8.8535465e-01 & -2.1532391e+01 \\
$y^3$   &   -7.6829624e+01 &  9.6863629e+01 \\
$x^4$   &  -1.7341620e-04 &  4.7139766e-03 \\
$x^3 y$ &   -6.6303349e-03 &  2.6373297e-01 \\
$x^2 y^2$ & 1.1436030e+00  & 1.3005160e+00 \\
$xy^3$  &  1.0748060e+00  & 2.0694130e+01\\
$y^4$   &  4.8756465e+01 & -4.4042922e+01 \\
$x^5$   & -1.6917964e-04 &  2.2350495e-04 \\
$x^4 y$ &-3.9532197e-03  &-5.4684219e-03 \\
$x^3 y^2$ &-1.3394365e-02 & -1.8843724e-01 \\
$x^2 y^3$ & -4.8332239e-01 & -1.0073151e+00 \\
$xy^4$  & -8.7359590e-01  &-7.7981035e+00\\
$y^5$   & -1.2443114e+01 &  6.0564361e+00 \\
$xdB1/dt$ &5.4772374e-01  & 4.3265057e-05 \\
$ydB2/dt$ & 1.6821393e-05 &  5.4531875e-01\\
\end{tabular}
\end{ruledtabular}
\end{table}

\begin{figure*}[htb]
\includegraphics[scale=0.20]{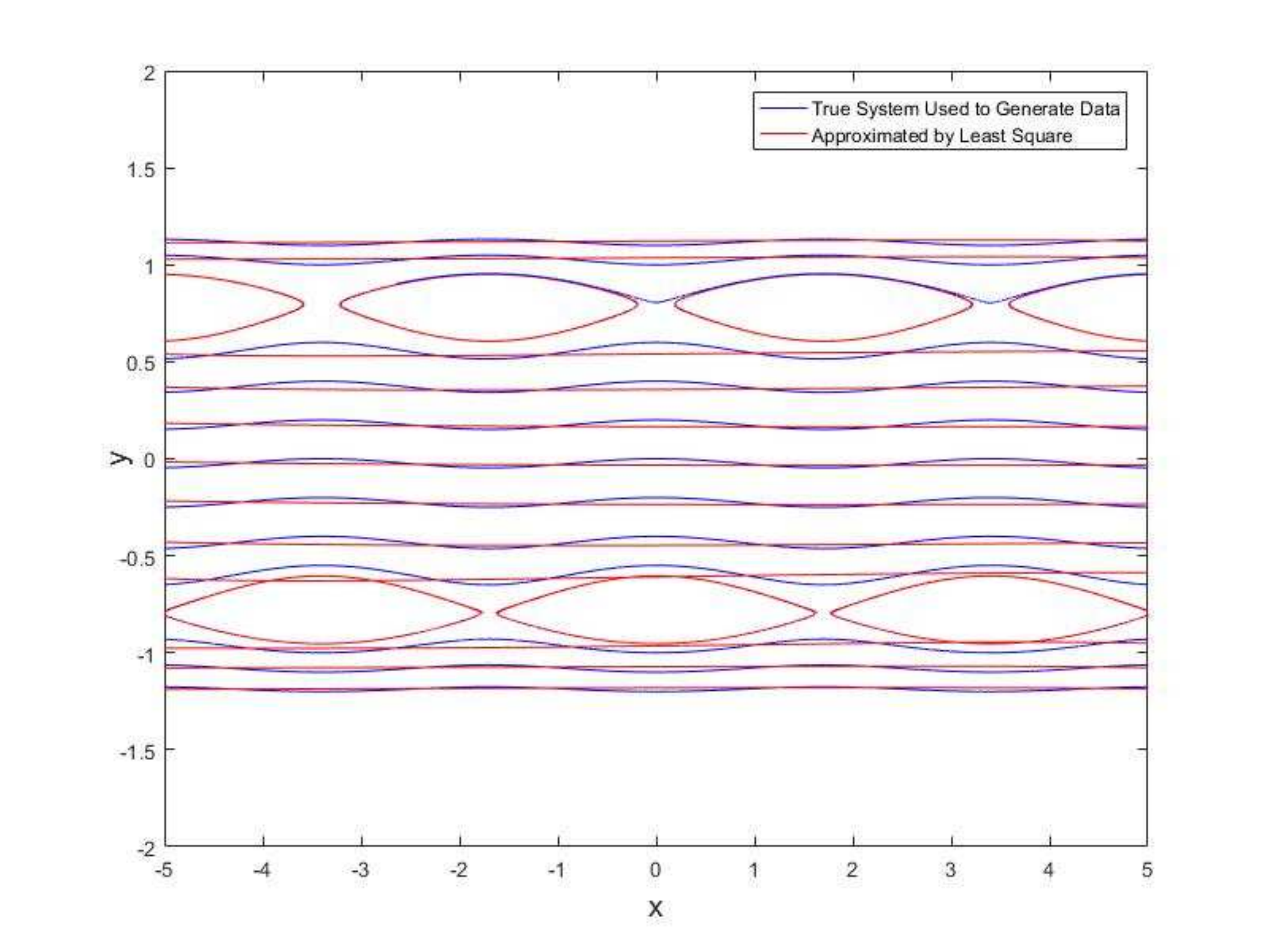}
\includegraphics[scale=0.20]{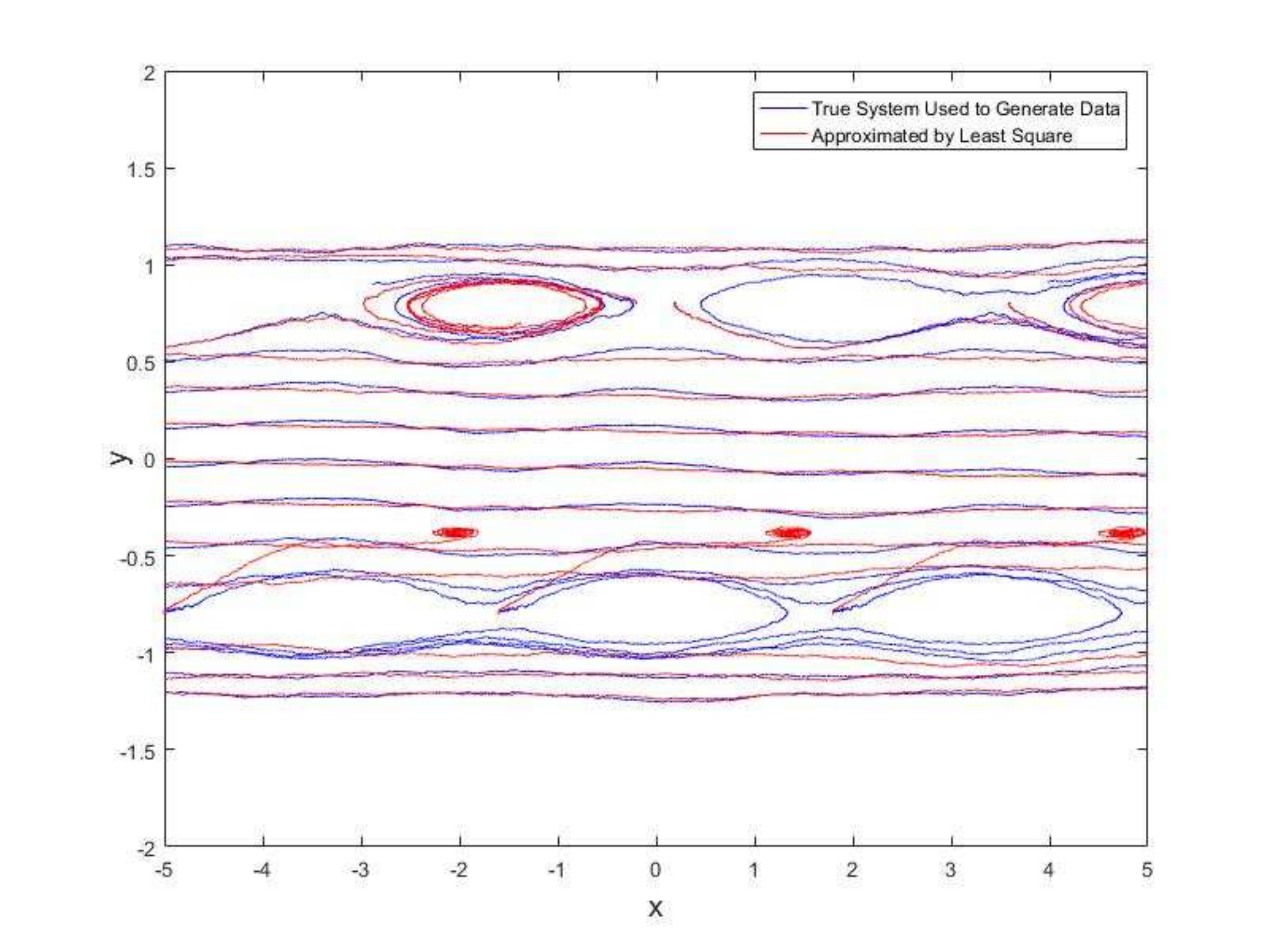}
\includegraphics[scale=0.20]{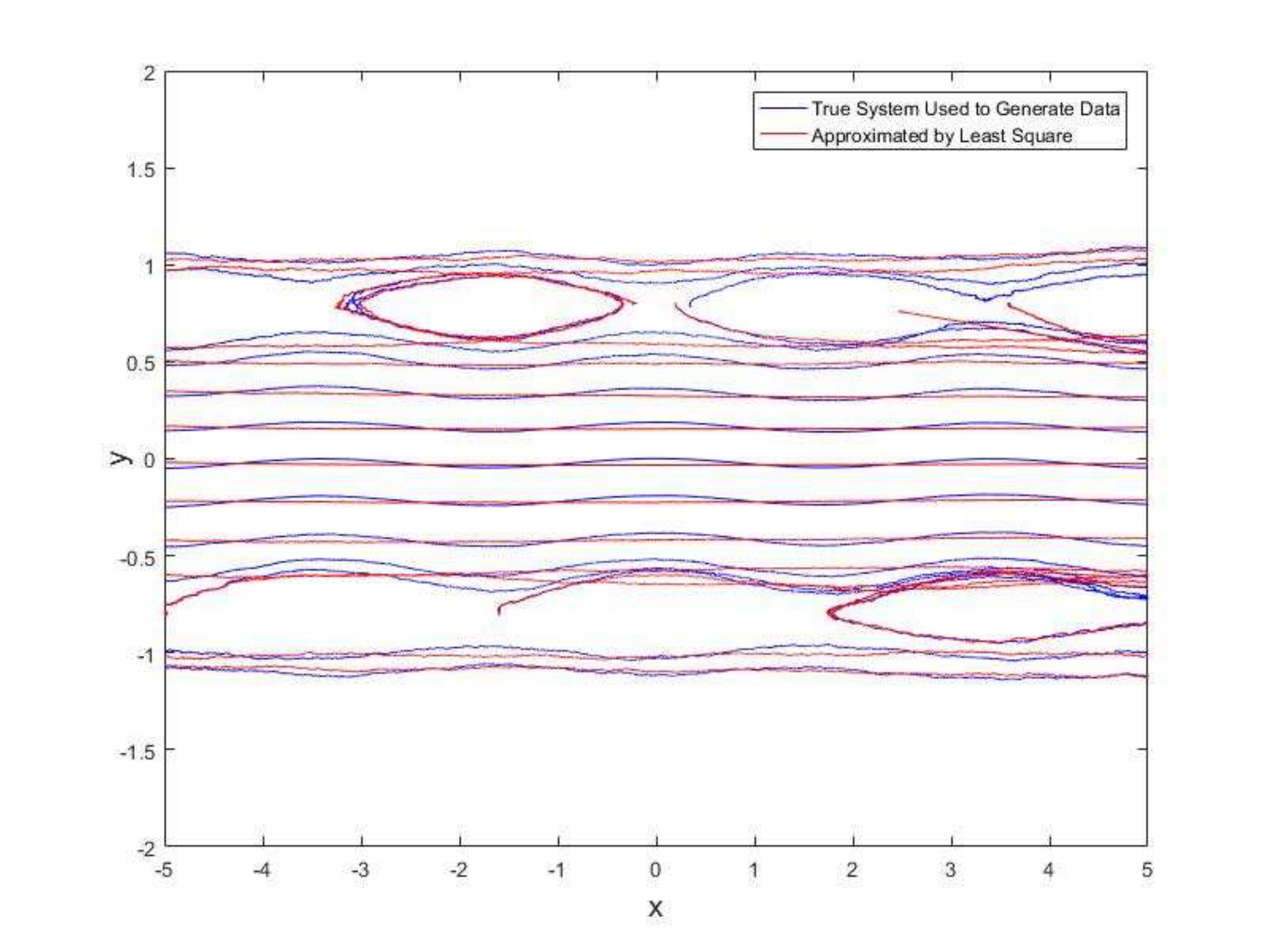}
\caption{Lagrangian fluid trajectories of the learned  meandering jet model  \eqref{learn1}- \eqref{learn2} with 20 initial positions: The fluid motion in the deterministic basic  model when noise is absent (left)
is steady, with no exchange of fluid between the eddies and the jet stream.
Under additive noise (middle) and multiplicative noise  (right),
there are   fluid exchanges between   eddies and the jet stream. }
\label{figure 2}
\end{figure*}

FIG. \ref{figure 2} shows  Lagrangian fluid trajectories of the learned  meandering jet model  \eqref{learn1}- \eqref{learn2}, without noise, with additive noise and with multiplicative noise.

\subsection{\label{sec:level2}Mean residence time }
\setcounter{equation}{0}
\renewcommand\theequation{3.\arabic{equation}}
In this subsection, we will consider  stochastic quasigeostrophic meandering  jet models,  \eqref{ad1}- \eqref{ad2} and \eqref{mu1}-\eqref{mu2},
to discover mean residence time by our   approach.
\begin{figure}[htb]
\center
\includegraphics[scale=0.25]{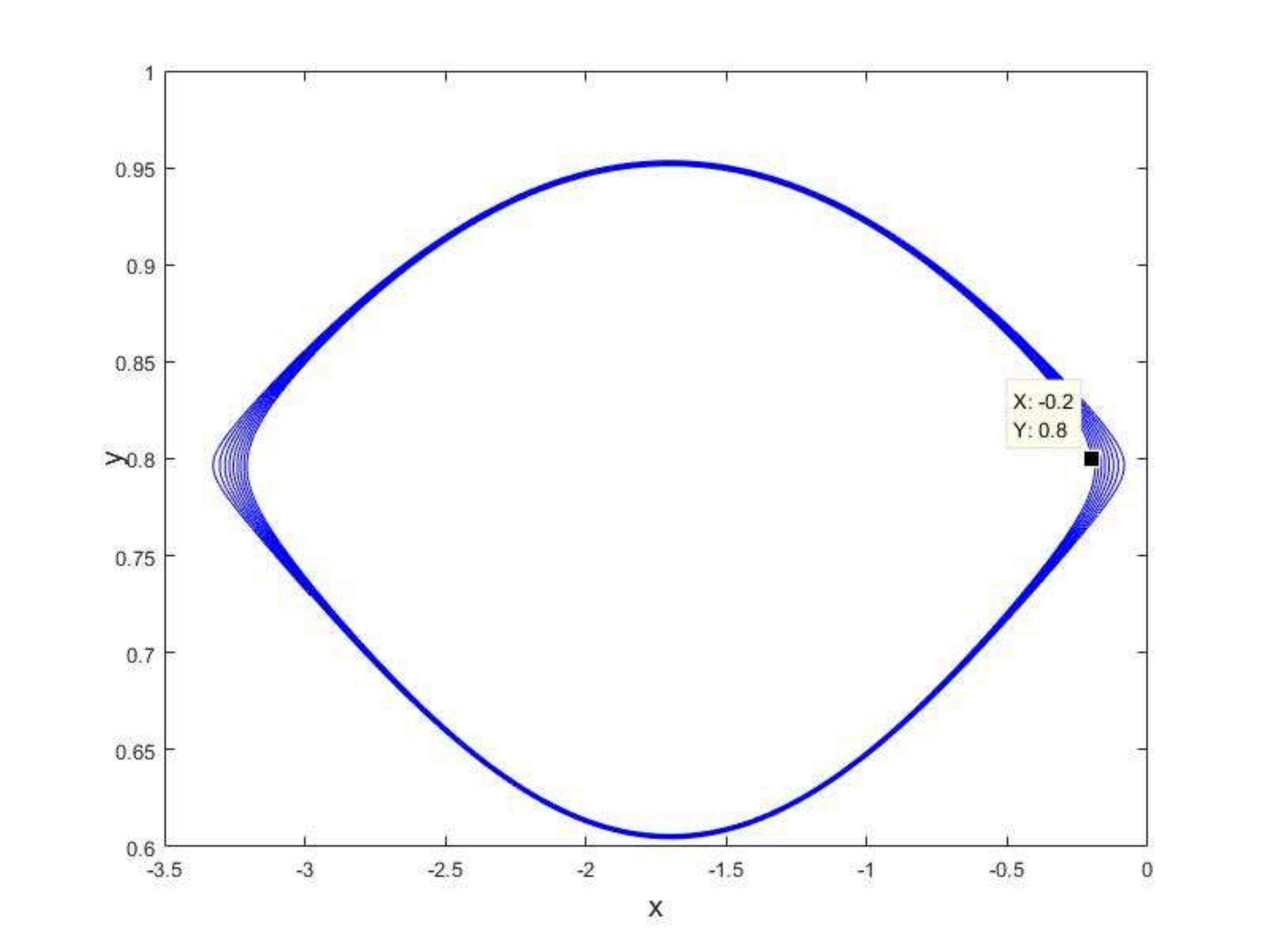}
\caption{An eddy: $\sigma=0$ .}
\label{figure 3}
\end{figure}
We take an eddy (FIG. \ref{figure 3}) as  the   bounded domain $D$ (with boundary $\partial D$ composed with
meander trough and crest).

The mean residence time $u(x,y)$ of the stochastic system \eqref{ad1}-\eqref{ad2} with additive noise
for a trajectory starting in the eddy, satisfies the following elliptic partial  differential equation as in \eqref{u1}-\eqref{u2}:
\begin{eqnarray}
\frac{1}{2}\sigma^2 \triangle u+(-\psi_{y})u_{x}+(\psi_{x})u_{y}=-1,\label{meanad1}\\
u|_{\partial D}=0.\label{meanad2}
\end{eqnarray}
The  mean residence time, as solution of  this  elliptic partial  differential equation, is denoted by $u_{True}$.
For the  corresponding learned    meandering jet model  \eqref{learn1}- \eqref{learn2}, with learned coefficients in Table I, we can also set up a similar  elliptic partial differential equation for learned mean residence time $u_{Learn}$.

Similarly, the mean residence time $u(x,y)$ of the stochastic system \eqref{mu1}-\eqref{mu2} with multiplicative  noise, for a trajectory starting in the eddy, satisfies the following elliptic partial  differential equation as in \eqref{u1}-\eqref{u2}:
\begin{eqnarray}
\frac{1}{2}\sigma^2 x^2 u_{xx}+ \frac{1}{2}\sigma^2 y^2 u_{yy}+(-\psi_{y})u_{x}+(\psi_{x})u_{y}=-1,\label{meanmu1}\\
u|_{\partial D}=0.\label{meanmu2}
\end{eqnarray}
The  mean residence time, as solution of  this  elliptic partial  differential equation, is denoted by $u_{True}$.
For the  corresponding learned    meandering jet model  \eqref{learn1}- \eqref{learn2}, with learned coefficients in Table II, we can also set up a similar elliptic partial differential equation for learned mean residence time $u_{Learn}$.



\begin{figure}[htb]
\center
\includegraphics[scale=0.15]{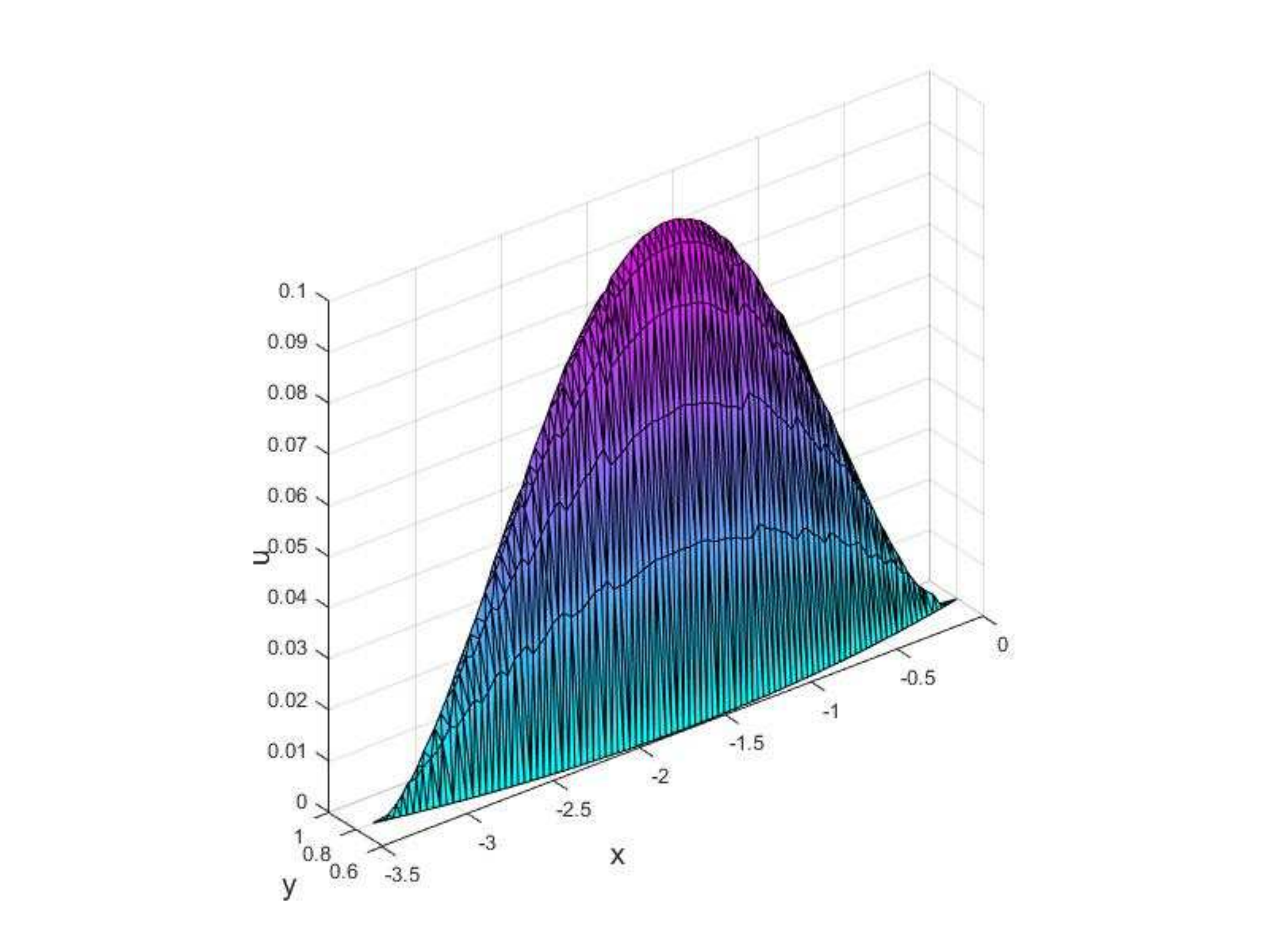}
\includegraphics[scale=0.15]{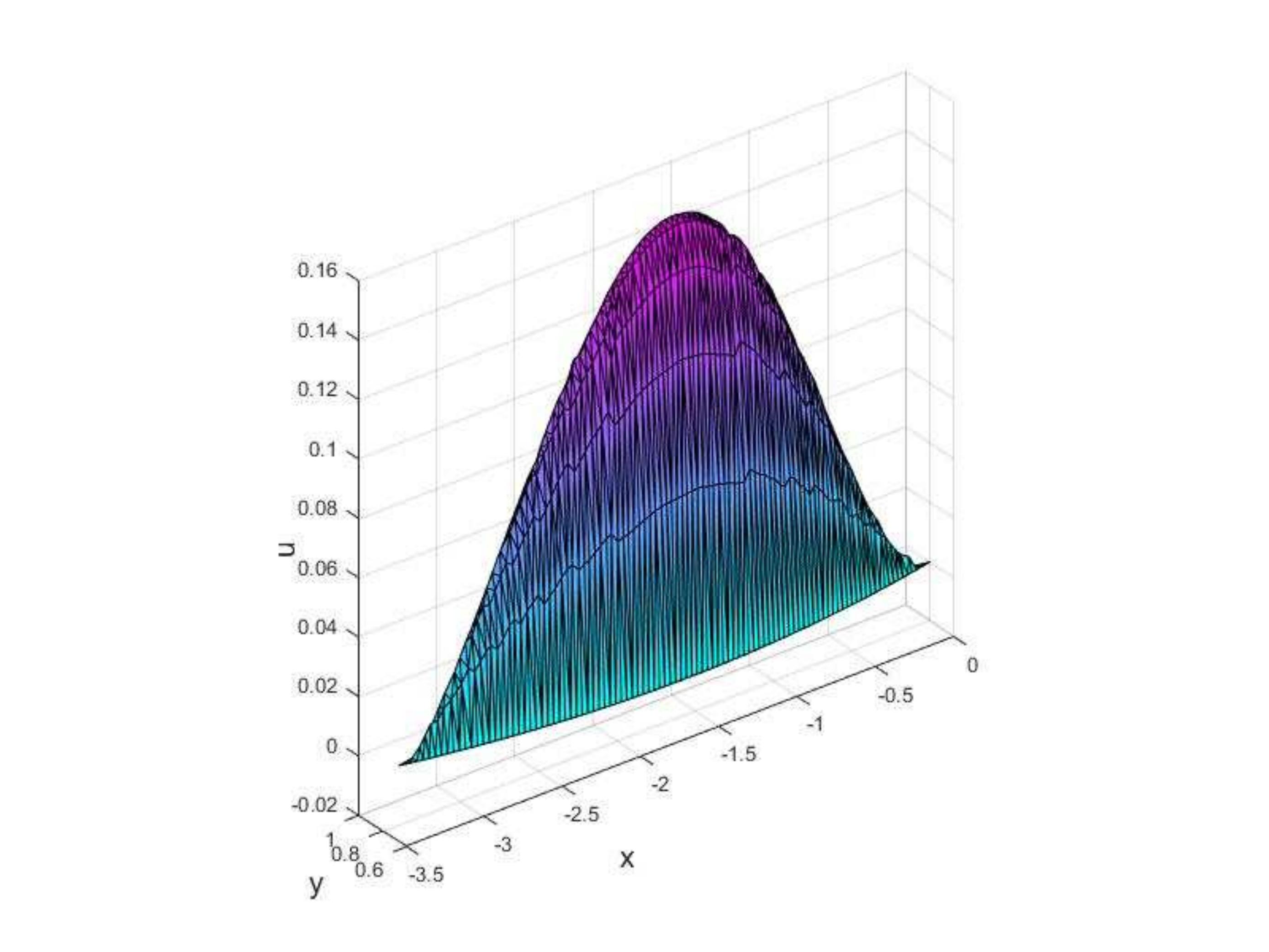}
\caption{Learned mean residence time $u_{Learn}$ for stochastic system \eqref{learn1}-\eqref{learn2} in an eddy.
Additive noise case is on the left,
 and multiplicative noise case is on the right.}
\label{figure 4}
\end{figure}

We use a finite element code  to solve the preceding  elliptic differential equations to get mean residence time
$u_{True}$ and  $u_{Learn}$, for both additive and multiplicative noise cases.   FIG.  \ref{figure 4} shows the mean residence time $u_{Learn}$, as  $u_{True}$ is barely distinguishable from $u_{Learn}$ when they are plotted together.  Thus we compare the error   between  the maximal values of the  mean residence times  $u_{True}$ and   $u_{Learn}$.
 In fact,  the error of  mean residence time between the learned system and the original system
 is $\mbox{error}(u)=\max_{(x,y)\in D}\|u_{Learn}-u_{True}\|= 2.2142\times 10^{-4}$ for additive noise case,
 and $\mbox{error}(u)=\max_{(x,y)\in D}\|u_{Learn}-u_{True}\|=1.3\times 10^{-3}$ for multiplicative noise case.

\subsection{\label{sec:level2} Escape probability}
\setcounter{equation}{0}
\renewcommand\theequation{4.\arabic{equation}}

We  again  take the eddy  (FIG. \ref{figure 3})  to be the   bounded domain $D$ (with boundary $\partial D$ composed with meander trough  (lower subboundary) and crest (upper subboundary)).
For the stochastic meandering jet system \eqref{ad1}-\eqref{ad2} with additive noise,  the escape probability $p(x,y)$ of a fluid particle,  starting at  $(x, y)$ in an eddy and escaping through  a  boundary component $\Gamma$ (either  the  trough or  crest),  satisfies the elliptic partial differential equation \eqref{p1}-\eqref{p3}:
\begin{eqnarray}
\frac{1}{2}\sigma^2 \triangle p+(-\psi_{y})p_{x}+(\psi_{x})p_{y}=0,\label{probad1}\\
p|_{\Gamma}=1,\label{probad2}\\
p|_{\partial D \setminus \Gamma}=0.\label{probad3}
\end{eqnarray}
For the  corresponding learned    meandering jet model  \eqref{learn1}- \eqref{learn2}, with learned coefficients in Table I, we can also set up a similar  elliptic partial differential equation for the  learned  escape probability  $p_{Learn}$.

Similarly,  for the stochastic meandering jet  system \eqref{mu1}-\eqref{mu2} with multiplicative noise, the escape probability $p(x,y)$ of   a fluid particle,  starting at  $(x, y)$ in an eddy and escaping through a  boundary component $\Gamma$ (either  the  trough or  crest),  satisfies the elliptic partial differential equation \eqref{p1}-\eqref{p3}:
\begin{eqnarray}
\frac{1}{2}\sigma^2 x^2 p_{xx}+ \frac{1}{2}\sigma^2 y^2 p_{yy}+(-\psi_{y})p_{x}+(\psi_{x})p_{y}=0,\label{probmu1}\\
p|_{\Gamma}=1,\label{probmu2}\\
p|_{\partial D \setminus \Gamma}=0.\label{probmu3}
\end{eqnarray}
For the  corresponding learned    meandering jet model  \eqref{learn1}- \eqref{learn2}, with learned coefficients in Table II, we can also set up a similar elliptic partial differential equation for the  learned   escape probability $p_{Learn}$.

\begin{figure}[htb]
\center
\includegraphics[scale=0.15]{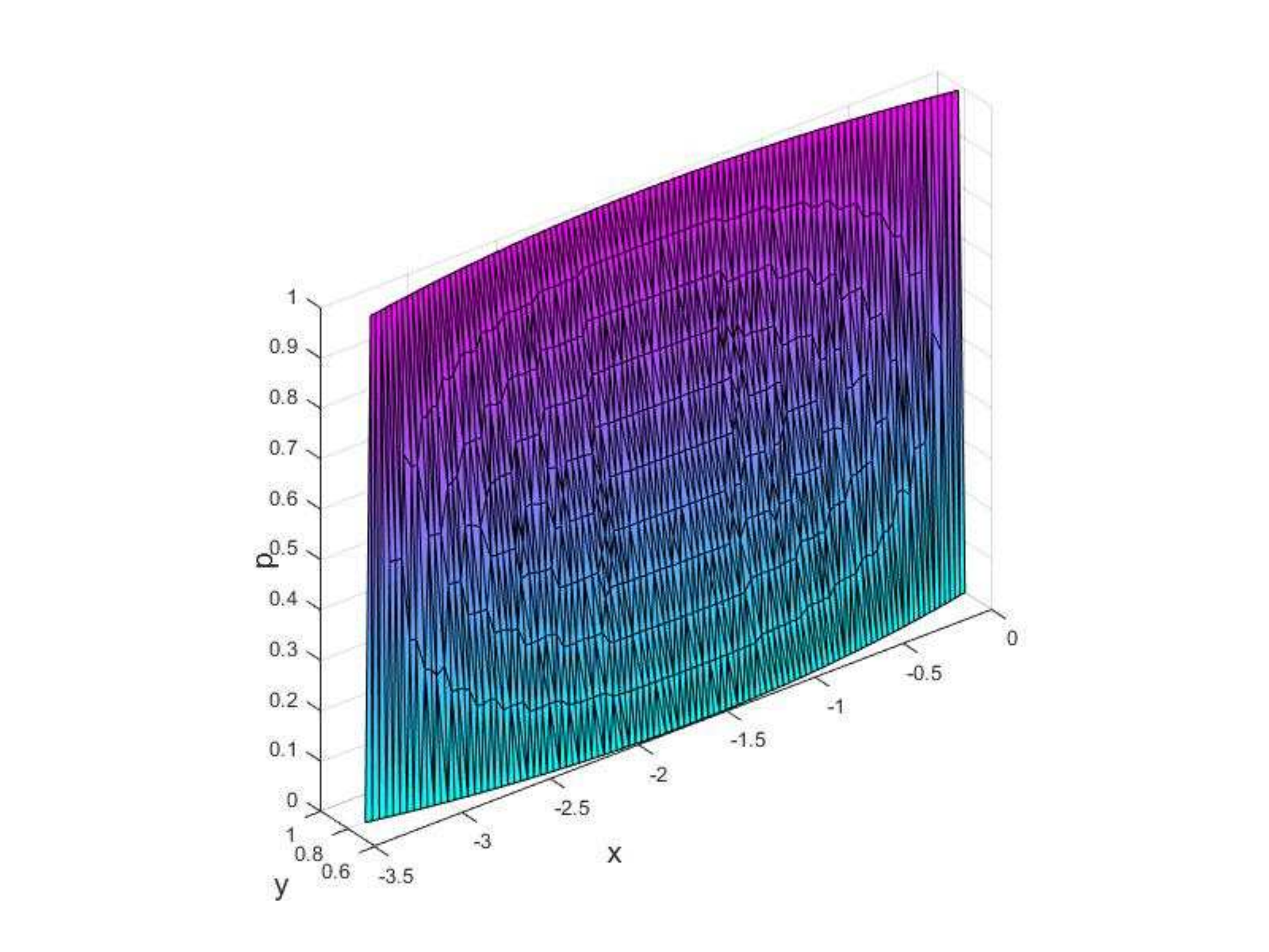}
\includegraphics[scale=0.15]{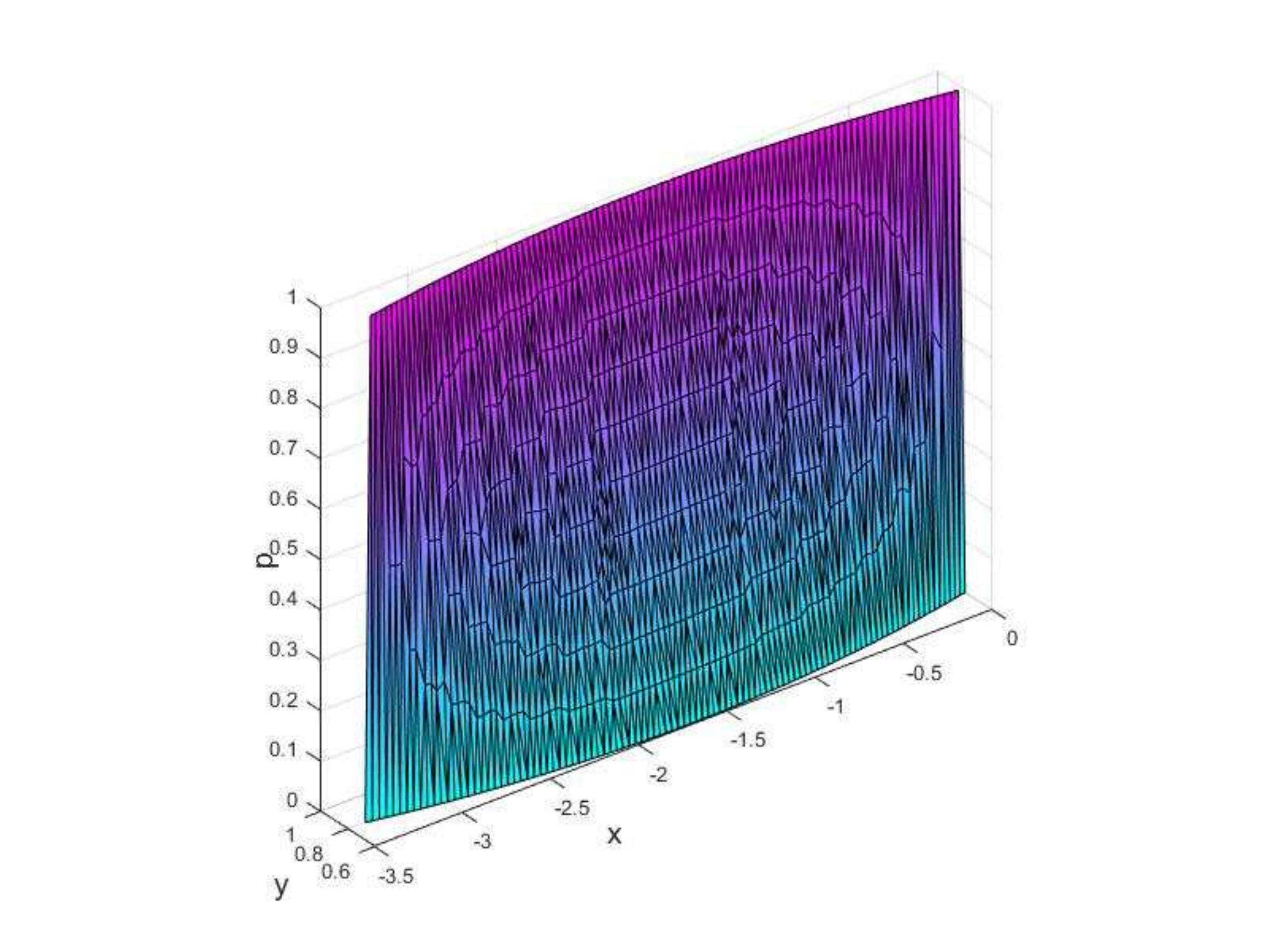}
\caption{Learned escape probability $P_{Learn}$ for stochastic system \eqref{learn1}-\eqref{learn2} in an eddy, exiting from upper subboundary.
 Additive noise case is on the left,
 and  multiplicative noise case is on the right.}
\label{figure 5}
\end{figure}
\begin{figure}[htb]
\center
\includegraphics[scale=0.15]{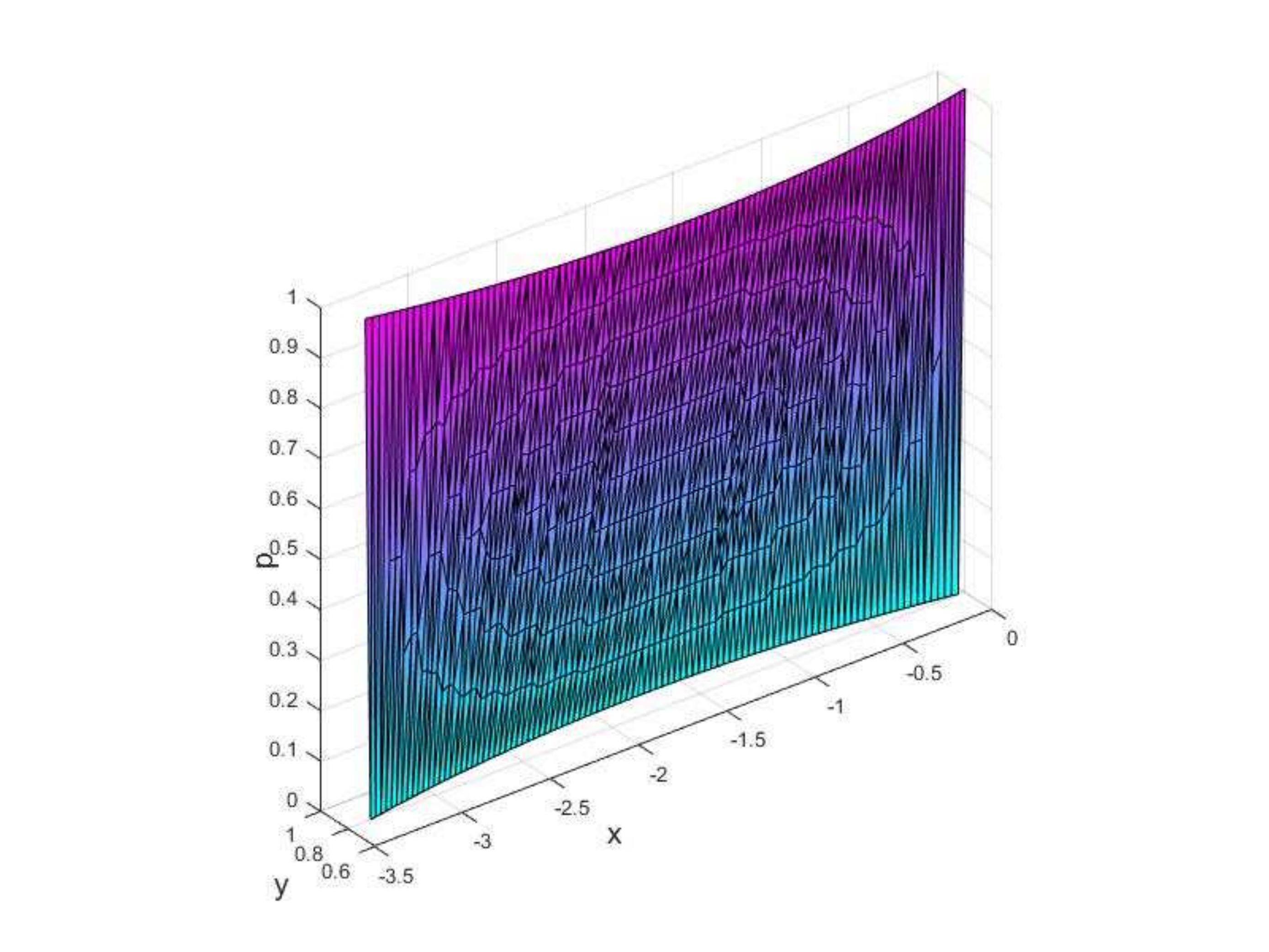}
\includegraphics[scale=0.15]{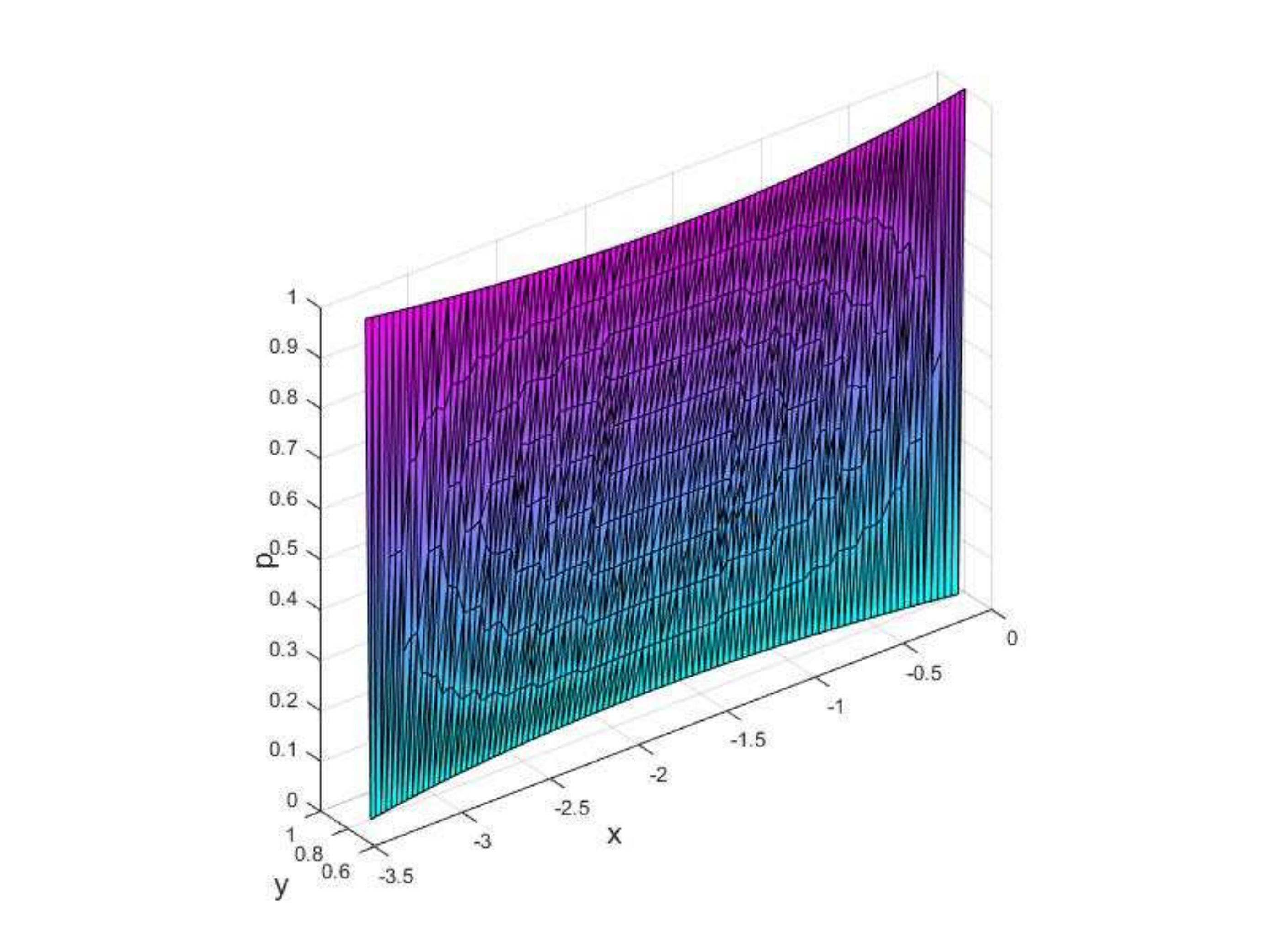}
\caption{Learned escape probability $P_{Learn}$ for stochastic system \eqref{learn1}-\eqref{learn2} in an eddy,
exiting from lower subboundary.
Additive noise case is on the left,
 and multiplicative noise case is on the right.}
\label{figure 6}
\end{figure}

FIGs. \ref{figure 5}-\ref{figure 6}  show the finite element solution for the  learned  escape probability  $p_{Learn}$ alone,  as  $p_{True}$ is barely distinguishable from $p_{Learn}$ when they are plotted together.  Thus we compare the error   between  the average  escape probability values.

 For fluid particles   initially uniformly distributed   in $ D $ (an eddy),
we     compute the average escape probability $P $ for   particles   leaving  $D$  through  the
upper or lower subboundary $\Gamma$,  given by $P = \frac{1}{D}\int\int_{D}p(x,y)dxdy$.
The error of average escape probability,  for fluid particles exiting fom $D$ through the upper subboundary,
between the learned system and  the original  system
is $\mbox{error}(P)=\|P_{Learn}-P_{True}\|=1.3949 \times 10^{-5}$ both for additive noise and multiplicative noise cases.
The error of the average escape probability,  for   particles   leaving  $D$  through the lower subboundary,
between the learned system and the original system is $\mbox{error}(P)=\|P_{Learn}-P_{True}\|=7.4133 \times 10^{-4}$
both for additive noise and multiplicative noise cases.

\section{\label{sec:level1} Learning three dimensional stochastic dynamical systems}
\setcounter{equation}{0}
\renewcommand\theequation{5.\arabic{equation}}

In this section, we will illustrate  our method in  two systems  with additive noise.
The first   system  is  for a three dimensional stochastic  linear   oscillator. The second
is  the more complex stochastic Lorenz system.
In these two examples, data from direct numerical simulations  (as `observation data')  are used
to discover   mean residence time and escape probability.

\subsection{\label{sec:level2} A stochastic linear system}

\begin{figure}[htb]
\center
\includegraphics[scale=0.15]{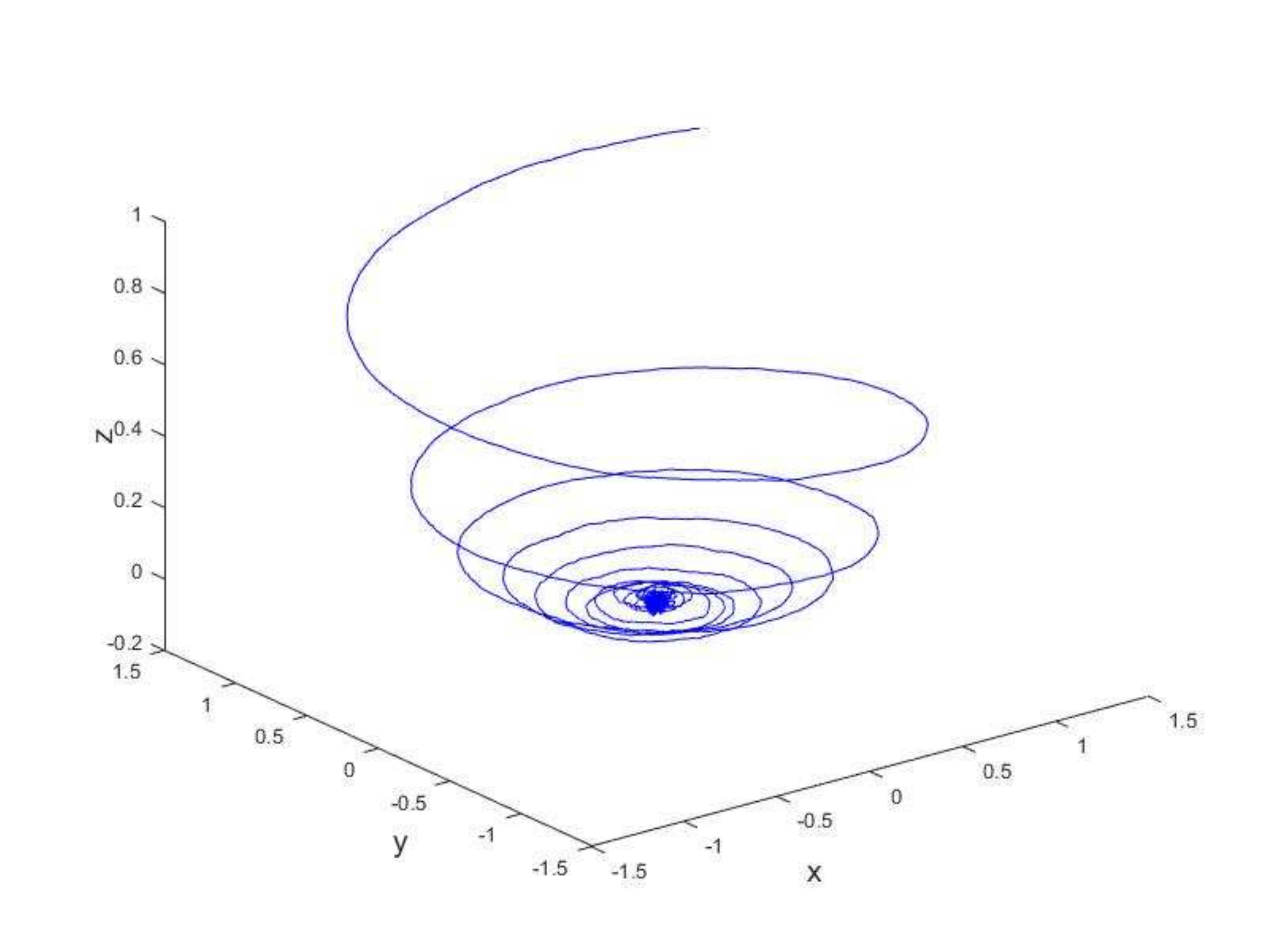}
\includegraphics[scale=0.15]{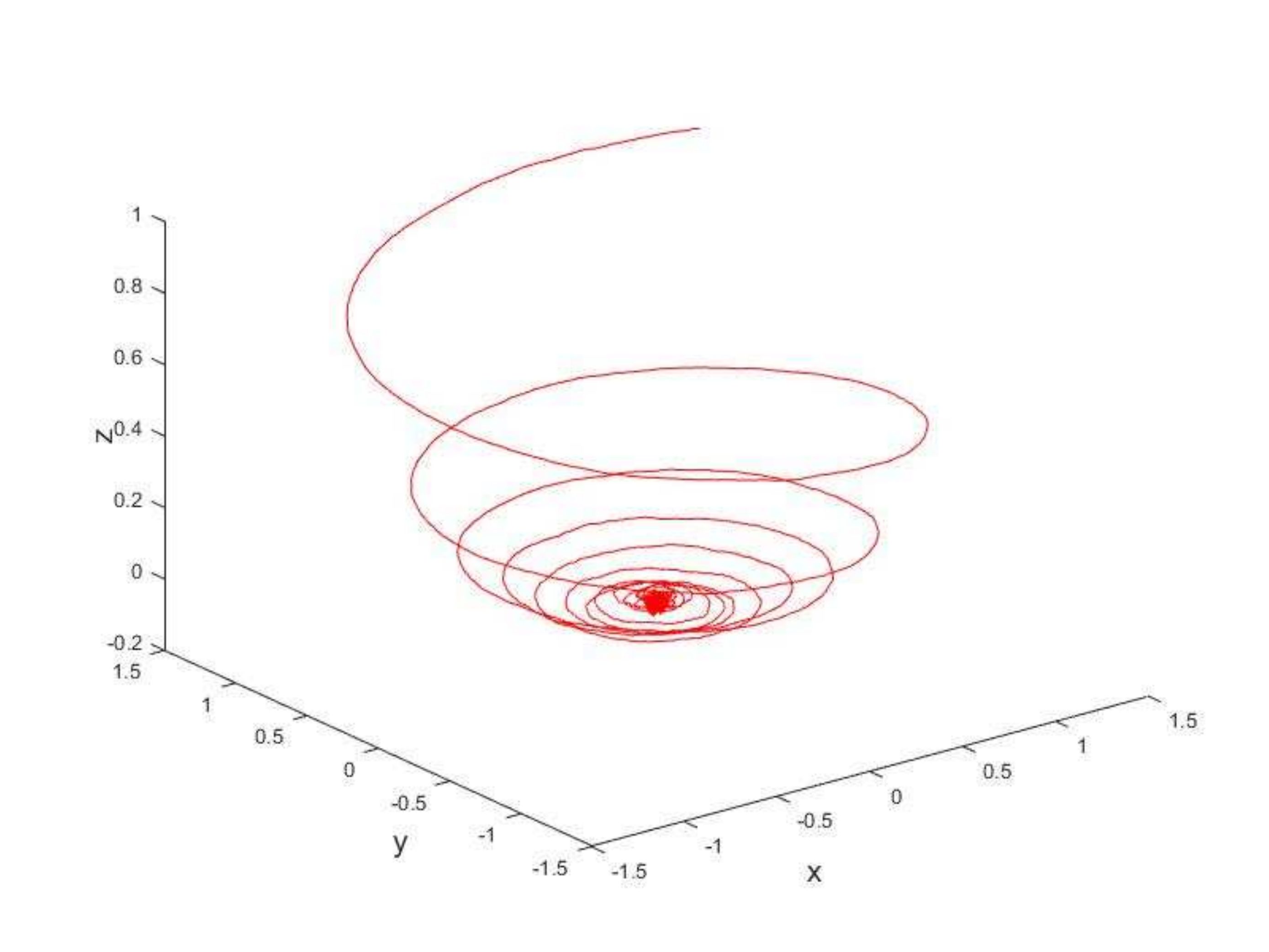}
\caption{Trajectory of  the three-dimensional stochastic linear system  starting at $(1,1,1)^{T}$:
The trajectory for the  original system in blue curve (left) and  for
the  learned system in red curve (right). }
\label{figure 7}
\end{figure}

We consider a linear stochastic system
\begin{eqnarray}\label{3Dlinear}
dx &=& (-0.1 x-2 y)dt + \sqrt{\epsilon} dB_1,\label{3Dlinear1}\\
dy &=& (2 x-0.1 y)dt + \sqrt{\epsilon} dB_2,\label{3Dlinear2}\\
dz &=& (-0.3 z)dt + \sqrt{\epsilon} dB_3\label{3Dlinear3},
\end{eqnarray}
where  $\epsilon$  is the  noise  intensity (taken to be $0.9$ here), and $B_1(t)$,  $B_2(t)$, $B_3(t)$ are three independent Brownian motions.
With  sample path data for this system, we try to  discover the  governing    equation,
\begin{eqnarray}
\dot{x} &=& \Theta \xi^1,  \label{3Dlearn1} \\
\dot{y} &=& \Theta \xi^2,\\
\dot{z} &=& \Theta \xi^3.  \label{3Dlearn3}
\end{eqnarray}
as in Brunton et al.\cite{Kutz} but we use a stochastic basis.  Here $\Theta$ is a set of  basis functions and $\Xi=[\xi^1,\xi^2,\xi^3]$ are coefficients (or weights).
As   in Section 2, we collect data $(x(t_i),y(t_i),z(t_i)),i=1,\cdots,20000$, from samplewise simulations of the original  system   \eqref{3Dlinear1}-\eqref{3Dlinear3},       and
construct a library $\Theta$ consisting of   polynomial   $\{1,x,y,z,\cdots, z^4\}$
and time derivatives of Brownian motions $\{dB_1/dt,dB_2/dt,dB_3/dt\}$.  We then solve  a  regression problem to determine the weights
$\xi^k=[\xi_1^k,\cdots,\xi_{38}^k],k=1,2,3.$  See Table \ref{tab:2} and FIG. \ref{figure 7}.

\begin{table}[h]
\caption{Identified coefficients for discovering a three-dimensional stochastic linear system
 with basis functions.  Sample paths of  \eqref{3Dlinear1}-\eqref{3Dlinear3} with initial position $(1,1,1)^{T}$   are numerically generated,    for $t\in [0,200]$ with stepsize $0.01$.}
\label{tab:2}
\begin{ruledtabular}
\begin{tabular}{cccc}
basis & $\dot{x}$ & $\dot{y}$& $\dot{z}$\\
\noalign{\smallskip}\hline\noalign{\smallskip}
$1$     &  -3.9018726e-06 & -2.5270553e-06 &  9.9621283e-06 \\
$x$     &  -9.9972593e-02 &  2.0000178e+00 & -6.9973977e-05 \\
$y$     &  -1.9999463e+00 & -9.9965198e-02 & -1.3719725e-04 \\
$z$     &   1.9721074e-04 &  1.2772392e-04 & -3.0050351e-01 \\
$x^2$   &   4.1955355e-05 &  2.7172467e-05 & -1.0711898e-04 \\
$xy$    &  -2.9428922e-06 & -1.9059698e-06 &  7.5136920e-06 \\
$xz$    &   1.2255057e-04 &  7.9370113e-05 & -3.1289194e-04 \\
$y^2$   &   2.4067948e-05 &  1.5587653e-05 & -6.1449464e-05 \\
$yz$    &   3.3875553e-04 &  2.1939567e-04 & -8.6489908e-04 \\
$z^2$   &   4.0214283e-03 &  2.6044858e-03 & -1.0267374e-02 \\
$x^3$   &  -1.1648871e-04 & -7.5444137e-05 &  2.9741500e-04 \\
$x^2 y$ &  -2.3616373e-04 & -1.5295190e-04 &  6.0296519e-04 \\
$x^2z$  &  -8.2888236e-04 & -5.3682726e-04 &  2.1162742e-03 \\
$xy^2$  &  -1.1776806e-04 & -7.6272713e-05 &  3.0068141e-04 \\
$xyz$   &   3.5179026e-05 &  2.2783764e-05 & -8.9817889e-05 \\
$xz^2$  &  -1.8287612e-03 & -1.1844007e-03 &  4.6691307e-03 \\
$y^3$   &  -2.3525676e-04 & -1.5236450e-04 &  6.0064954e-04 \\
$y^2z$  &  -1.0099170e-03 & -6.5407468e-04 &  2.5784856e-03 \\
$yz^2$  &  -6.2014715e-03 & -4.0163950e-03 &  1.5833386e-02 \\
$z^3$   &  -1.2428487e-02 & -8.0493336e-03 &  3.1731989e-02 \\
$x^4$   &  -9.3740853e-05 & -6.0711444e-05 &  2.3933595e-04 \\
$x^3y$  &  -2.1022902e-06 & -1.3615523e-06 &  5.3674957e-06 \\
$x^3z$  &   7.5057766e-04 &  4.8611307e-04 & -1.9163493e-03 \\
$x^2y^2$&  -1.2711976e-04 & -8.2329355e-05 &  3.2455783e-04 \\
$x^2yz$ &   1.8010661e-03 &  1.1664639e-03 & -4.5984206e-03 \\
$x^2z^2$&   2.2180393e-03 &  1.4365174e-03 & -5.6630223e-03 \\
$xy^3$  &  -8.7412104e-06 & -5.6612618e-06 &  2.2317761e-05 \\
$xy^2z$ &   7.6527222e-04 &  4.9563003e-04 & -1.9538669e-03 \\
$xyz^2$ &  -3.0539398e-05 & -1.9778900e-05 &  7.7972151e-05 \\
$xz^3$  &   4.8388992e-05 &  3.1339224e-05 & -1.2354512e-04 \\
$y^4$   &  -3.2625868e-05 & -2.1130207e-05 &  8.3299255e-05 \\
$y^3z$  &   1.7990625e-03 &  1.1651663e-03 & -4.5933050e-03 \\
$y^2z^2$&   2.2767938e-03 &  1.4745699e-03 & -5.8130323e-03 \\
$yz^3$  &   3.0332728e-03 &  1.9645050e-03 & -7.7444488e-03 \\
$z^4$   &   5.7423101e-03 &  3.7190182e-03 & -1.4661070e-02 \\
$dB1/dt$&   9.4868476e-01 &  9.4448784e-07 & -3.7233490e-06 \\
$dB2/dt$&  -1.8088277e-05 &  9.4867158e-01 &  4.6182372e-05 \\
$dB2/dt$&   4.2468852e-05 &  2.7505034e-05 &  9.4857487e-01 \\
\end{tabular}
\end{ruledtabular}
\end{table}

We take cuboid $D = [-2,2] \times [-2,2] \times [0,1]$,
 with boundary $\partial D$,  containing a subboundary $\Gamma$ to be the
surface $z=1$ (top boundary) or the surface  $z=0$ (bottom boundary).

For the original system  \eqref{3Dlinear1}-\eqref{3Dlinear3},  the mean residence time $u$ (for solutions with initial points in $D$)  and escape probability $p$ (for solutions exiting through a subboundary $\Gamma$)  satisfy the elliptic partial
differential equations \eqref{u1}-\eqref{u2} and \eqref{p1}-\eqref{p3}, respectively. Hence,
\begin{eqnarray*}
\frac{1}{2}\epsilon \triangle u+(-0.1 x-2 y)u_{x}+(2 x-0.1 y)u_{y}+(-0.3 z)u_{z} &=& -1,\\
u|_{\partial D} &=& 0.
\end{eqnarray*}
\begin{eqnarray*}
\frac{1}{2}\epsilon \triangle p+(-0.1 x-2 y)p_{x}+(2 x-0.1 y)p_{y}+(-0.3 z)p_{z} &=& 0,\\
p|_{\Gamma} &=& 1,\\
p|_{\partial D \setminus \Gamma} &=& 0.
\end{eqnarray*}
These equations can also be solved by a finite element method to get the mean residence time $u_{True}$,  and the escape probability $p_{True}$ and the average escape probability $P_{True}$.

For the learned model  \eqref{3Dlearn1}-\eqref{3Dlearn3}, we also have  similar  elliptic partial differential equations for the mean residence time $u_{Learn}$,    the escape probability $p_{Learn}$ and the average escape probability $P_{Learn}$.



\begin{figure}[htb]
\center
\includegraphics[scale=0.3]{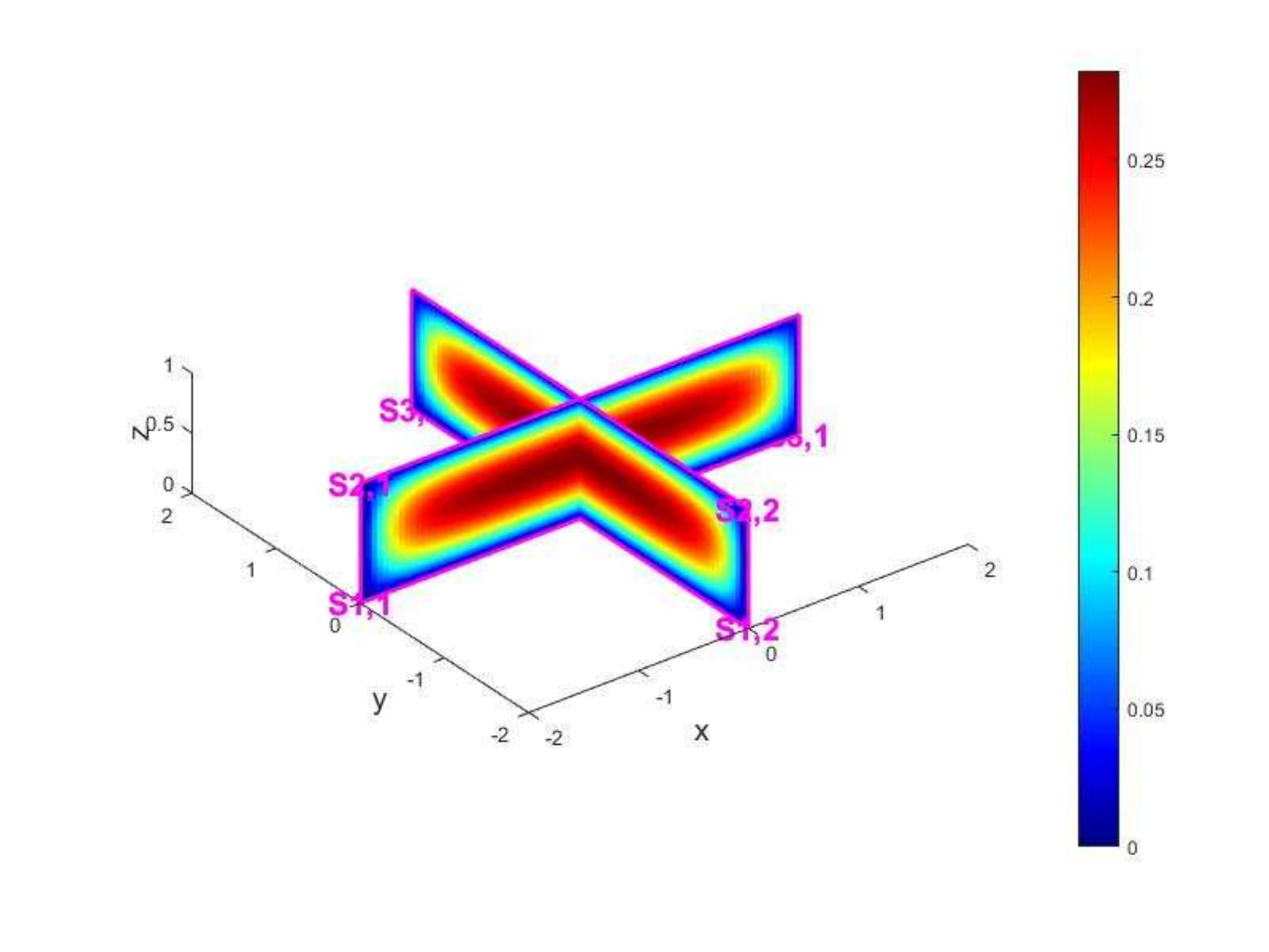}
\caption{Learned mean residence time $u_{Learn}$ for   the learned   linear stochastic system \eqref{3Dlearn1}-\eqref{3Dlearn3}, as viewed in two slices. }
\label{figure 8}
\end{figure}

The mean residence time $u(x,y,z)$ of the learned model is shown in
  FIG. \ref{figure 8}.  It is barely distinguishable from  the mean residence time for the original  system and we thus do not show the latter. Instead,  we compute the error (using 20000  uniformly  distributed  points in $D$)  between   the maximal values of   mean residence time for the learned and original systems:
  $\mbox{error}(u)=\max_{(x,y,z)\in D}\|u_{Learn}-u_{True}\|=3.7495\times 10^{-4}$.

\begin{figure}[htb]
\center
\includegraphics[scale=0.15]{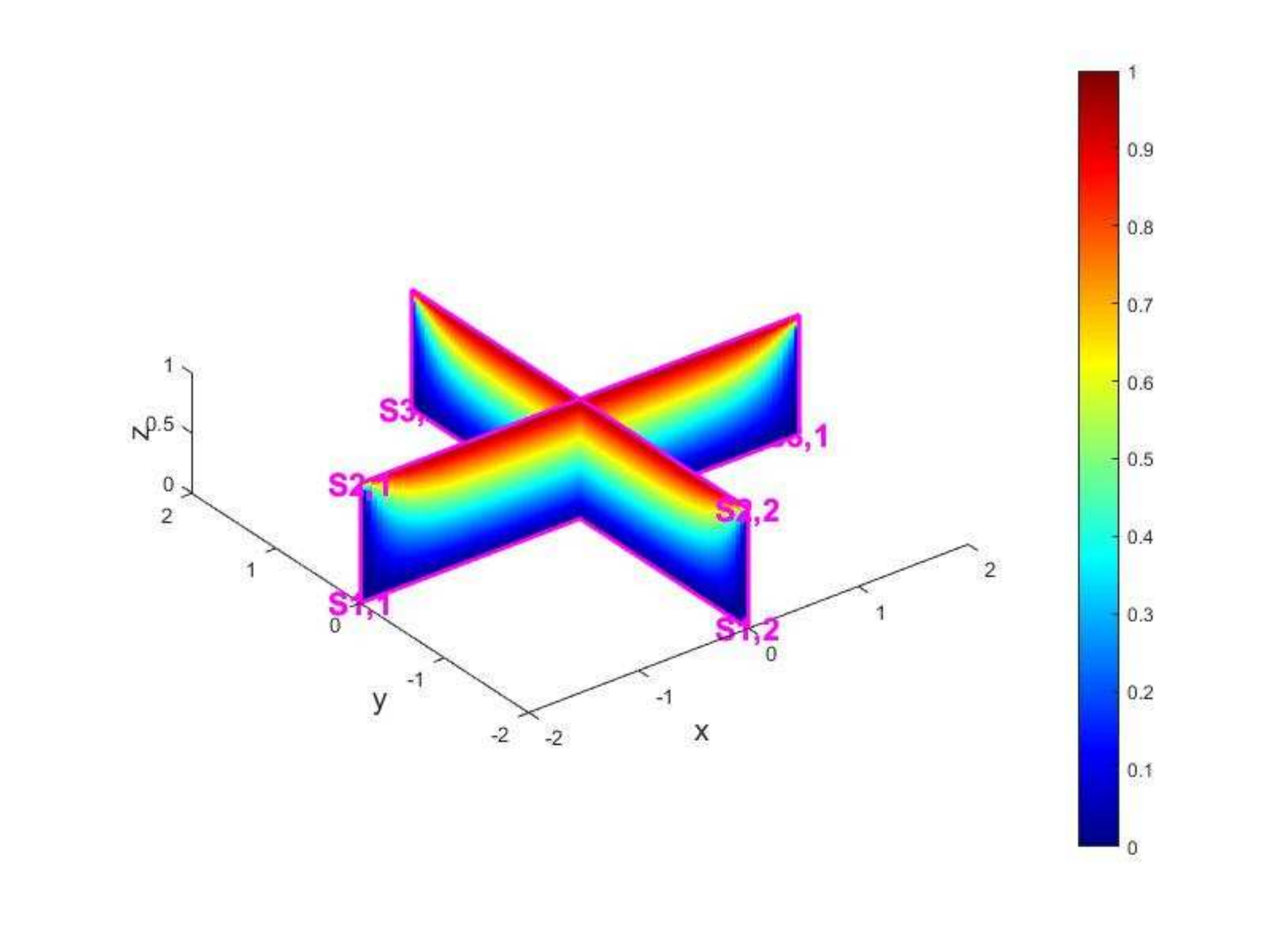}
\includegraphics[scale=0.15]{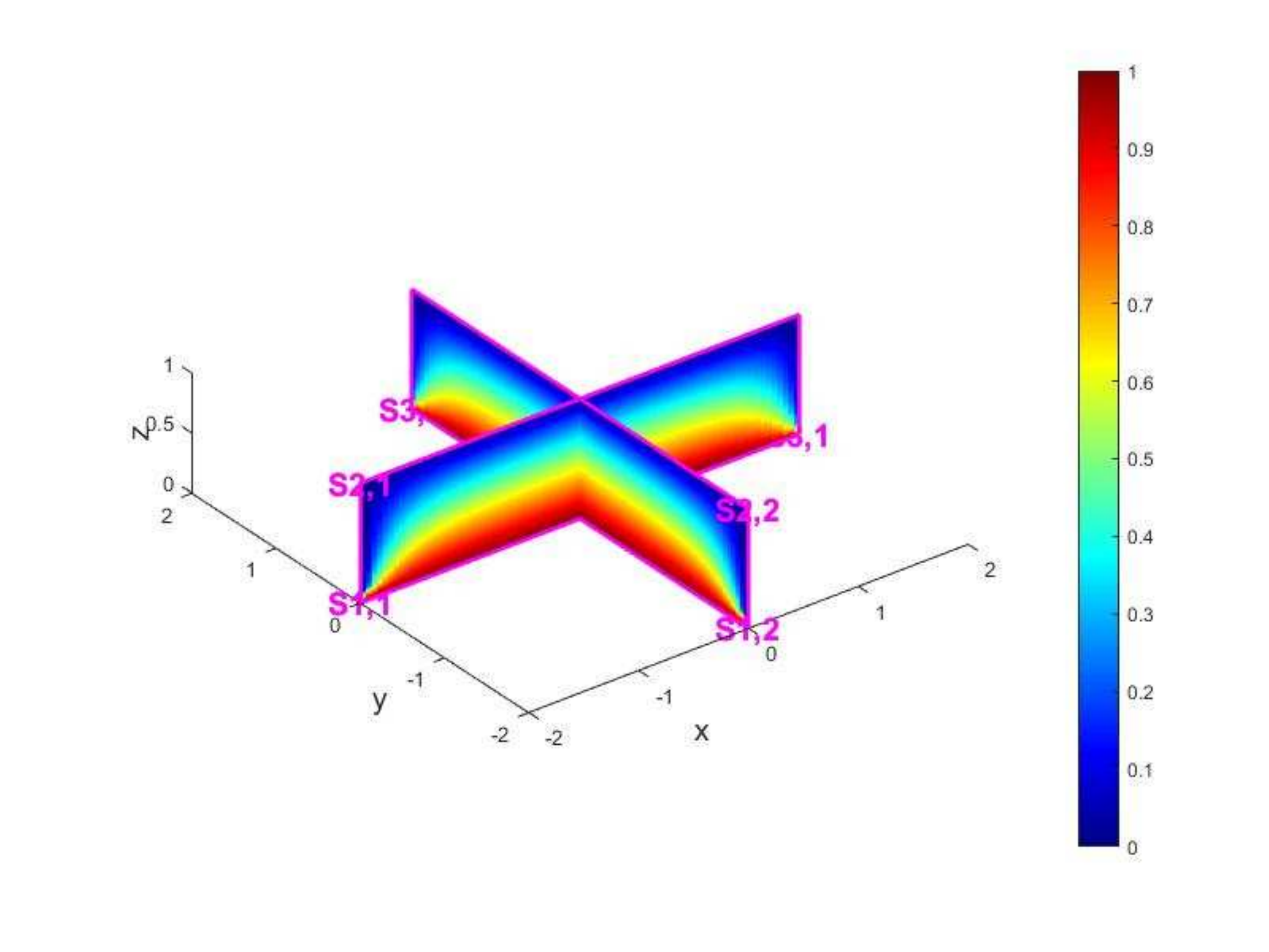}
\caption{Learned escape probability  for  the learned linear stochastic system  \eqref{3Dlearn1}-\eqref{3Dlearn3}:
 escaping from top surface (left) and bottom surface (right).  Two slices are shown for each subfigure.  }
\label{figure 9}
\end{figure}

The   escape probability    $p(x,y,z)$ of the learned system is shown in
  FIG. \ref{figure 9}.  It is barely distinguishable from  the  escape probability
  for the original system and we thus do not show the latter. Instead,  we calculate  the  error  between the  average escape probability values  of the
  learned system and original system is  $\mbox{error}(P)=\|P_{Learn}-P_{True}\|=9.9478 \times 10^{-7}$  (escaping from the top boundary)
and $\mbox{error}(P)=\|P_{Learn}-P_{True}\|=1.0783 \times 10^{-4}$ (escaping from the  bottom boundary).

\subsection{\label{sec:level2}Lorenz system with random noises}
We consider a stochastic  Lorenz  system \cite{Kutz}:
\begin{eqnarray}
dx &=& \sigma(y-x)dt + \sqrt{\epsilon} dB_1,    \label{Lorenz1}\\
dy &=& (\rho x-xz-y)dt + \sqrt{\epsilon} dB_2,   \label{Lorenz2}\\
dz &=& (xy-\beta z)dt + \sqrt{\epsilon} dB_3,     \label{Lorenz3}
\end{eqnarray}
where $\sigma=10$,  $\beta=\frac{8}{3}$,
$\rho=28$ are standard parameters,  and $B_1(t)$,  $B_2(t)$,  $B_3(t)$ are independent scalar Brownian motions.
We take $\epsilon=0.9$ in the  following computations.

As     before,  a  learned model is identified  in the following form
\begin{eqnarray}
\dot{x} &=& \Theta \xi^1,      \label{Lorenz4} \\
\dot{y} &=& \Theta \xi^2,         \label{Lorenz5}      \\
\dot{z} &=& \Theta \xi^3,      \label{Lorenz6}
\end{eqnarray}
where $\Theta$ is a set of  basis function consisting of polynomials   in
$(x,y,z)$ up to fourth order and time derivatives of three independent Brownian motions.
Then we  solve a regression problem to determine the weights $\Xi=[\xi^1,\xi^2,\xi^3]$; see Table \ref{tab:3}.

\begin{figure}[htb]
\center
\includegraphics[scale=0.15]{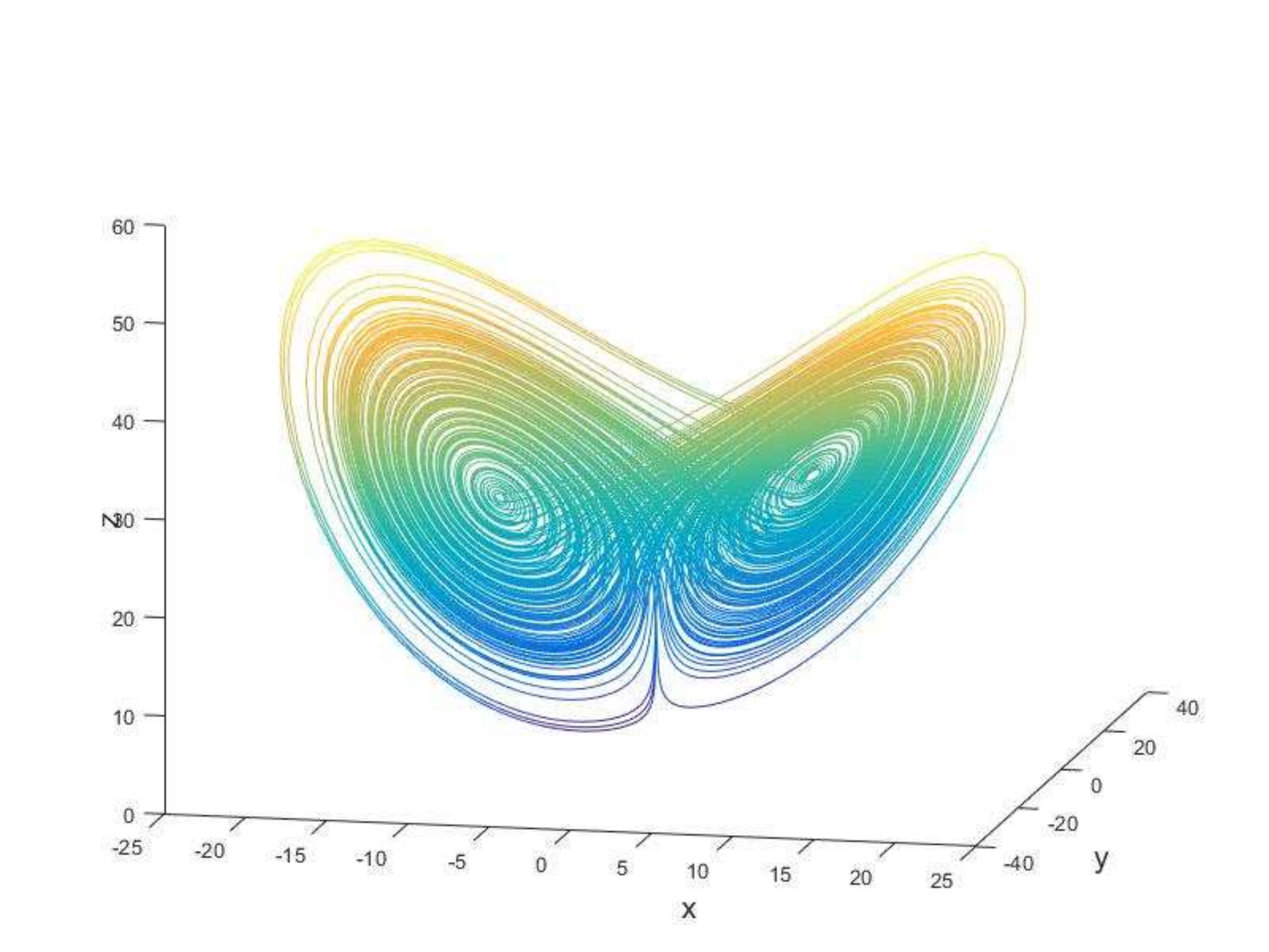}
\includegraphics[scale=0.15]{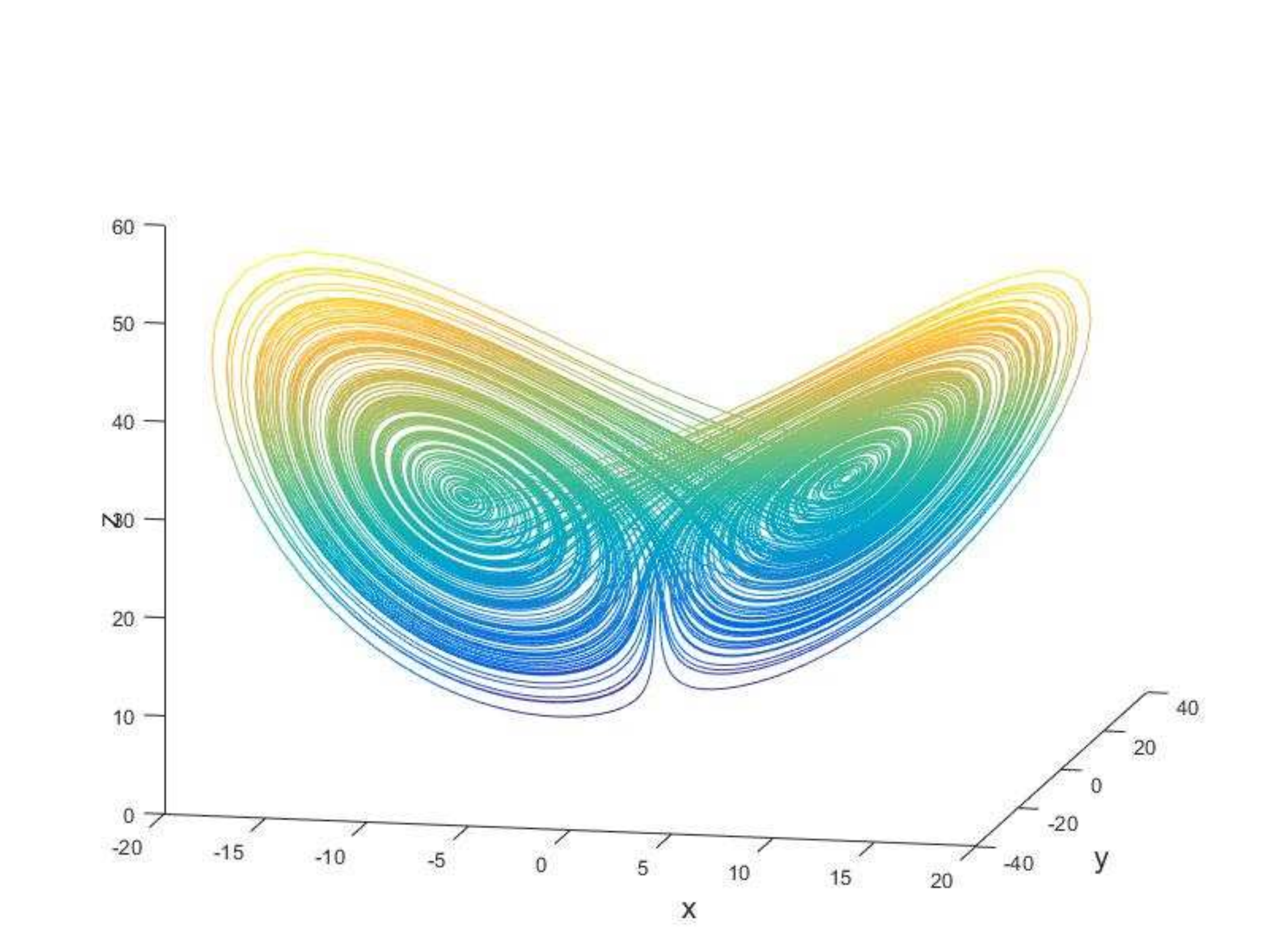}
\caption{One   trajectory  starting at $(x,y,z)^T=(-8,7,27)^T$:  for    the original
   Lorenz system (left), and for the
learned   Lorenz system  with basis in  Table \ref{tab:3} (right).}
\label{figure 10}
\end{figure}

\begin{table}[h]
\caption{ Identified Lorenz system
 with basis functions.   The sample paths  data are numerically generated
 by solving \eqref{Lorenz1}-\eqref{Lorenz3} with  initial position $(-8,7,27)^{T}$, for $t\in [0,200]$ with stepsize $0.01$.}
\label{tab:3}
\begin{ruledtabular}
\begin{tabular}{cccc}
basis & $\dot{x}$ & $\dot{y}$& $\dot{z}$\\
\noalign{\smallskip}\hline\noalign{\smallskip}
$1$     & 7.5197692e-03  & 3.3528379e-03 & -1.3696228e-02\\
$x$     & -1.0008895e+01 &  2.7996034e+01 & 1.6200350e-02\\
$y$     & 1.0004690e+01  &-9.9790907e-01  &-8.5413967e-03\\
$z$     &-2.7045815e-03 & -1.2058912e-03 & -2.6617406e+00\\
$x^2$   &-4.5633362e-04 & -2.0346537e-04  & 8.3114908e-04\\
$xy$    &-6.6325912e-04 & -2.9572721e-04  & 1.0012080e+00\\
$xz$   & 1.0139842e-03 & -9.9954790e-01 & -1.8468331e-03\\
$y^2$  &  4.1804010e-04 &  1.8639145e-04 & -7.6140269e-04\\
$yz$   & -5.5344389e-04 & -2.4676391e-04 &  1.0080221e-03\\
$z^2$  & 3.0860393e-04  & 1.3759717e-04 & -5.6207972e-04\\
$x^3$   & 3.8011238e-05  & 1.6948062e-05 & -6.9232255e-05\\
$x^2 y$ &-5.6250702e-05 & -2.5080489e-05  & 1.0245294e-04\\
$x^2z$ &4.3954760e-05  & 1.9598100e-05 & -8.0057566e-05\\
$xy^2$ & 1.1078797e-05  & 4.9397009e-06 & -2.0178509e-05\\
$xyz$  &1.6070760e-05   &7.1654662e-06 & -2.9270684e-05\\
$xz^2$  &-3.4384030e-05  &-1.5330800e-05 &  6.2625795e-05\\
$y^3$   &4.3084220e-06  & 1.9209952e-06  &-7.8471997e-06\\
$y^2z$ & -1.1913783e-05  &-5.3119961e-06  & 2.1699322e-05\\
$yz^2$  &2.0432252e-05   &9.1101240e-06  &-3.7214544e-05\\
$z^3$   & -1.3645677e-05 &-6.0841951e-06  & 2.4853729e-05\\
$x^4$   &2.4710563e-07  & 1.1017693e-07  &-4.5006900e-07\\
$x^3y$  &-1.5318615e-06 & -6.8301077e-07 &  2.7900756e-06\\
$x^3z$  &-6.0616355e-07  &-2.7027001e-07  & 1.1040438e-06\\
$x^2y^2$&  1.9963786e-06  & 8.9012488e-07 & -3.6361297e-06\\
$x^2yz$ &  1.1588151e-06 &  5.1668065e-07 & -2.1106228e-06\\
$x^2z^2$& -1.0342002e-06 & -4.6111863e-07  & 1.8836538e-06\\
$xy^3$  &  -7.4783630e-07 & -3.3343761e-07  & 1.3620812e-06\\
$xy^2z$ &  -2.2736871e-07 & -1.0137684e-07 &  4.1412091e-07\\
$xyz^2$ &2.5757823e-07 &  1.1484635e-07 & -4.6914341e-07\\
$xz^3$  & 3.3732532e-07 &  1.5040317e-07 & -6.1439179e-07\\
$y^4$   &-2.9520887e-09 & -1.3162472e-09  & 5.3768245e-09\\
$y^3z$  & -1.1011727e-07 & -4.9097963e-08 &  2.0056350e-07\\
$y^2z^2$&-6.2710697e-08 & -2.7960805e-08  & 1.1421893e-07\\
$yz^3$  &-2.3846556e-07 & -1.0632459e-07  & 4.3433230e-07\\
$z^4$   &2.0349814e-07 &  9.0733672e-08 & -3.7064394e-07\\
$dB1/dt$& 9.4581495e-01 & -1.2789082e-03  & 5.2242963e-03\\
$dB2/dt$&-3.1279105e-03 &  9.4728866e-01 &  5.6970599e-03\\
$dB3/dt$&  3.4942944e-03 &  1.5580003e-03  & 9.4231892e-01\\
\end{tabular}
\end{ruledtabular}
\end{table}

A single trajectory of this stochastic system with initial condition $(x,y,z)^T=(-8,7,27)^T$
is shown  in FIG. \ref{figure 10} (left), together with the same trajectory captured by the learned model in FIG. \ref{figure 10} (right).

For the stochastic  Lorenz system, we take a cuboid $D$ to be either
$D_1=[-9,-8]\times[-9,-8]\times[27,28]$  containing the  left saddle as left residence region,
or      $D_2=[8,9]\times[8,9]\times[27,28]$ containing the  right saddle as the right residence region.
Then mean residence time and escape probability satisfy the following elliptic partial
differential equations, respectively:
\begin{eqnarray*}
\frac{1}{2}\varepsilon \triangle u+(\sigma(y-x))u_{x}+(\rho x-xz-y)u_{y}+(xy-\beta z)u_{z} &=& -1,\\
u|_{\partial D_{i}} = 0,~i=1,2.
\end{eqnarray*}
\begin{eqnarray*}
\frac{1}{2}\varepsilon \triangle p+(\sigma(y-x))p_{x}+(\rho x-xz-y)p_{y}+(xy-\beta z)p_{z} &=& 0,\\
p|_{\Gamma} &=& 1,\\
p|_{\partial D_i \setminus \Gamma} = 0,~i=1,2.
\end{eqnarray*}
The finite element numerical  solutions to these partial differential  equations  are denoted by $u_{True}$ and $p_{True}$.

For the learned   Lorenz system    \eqref{Lorenz4}-\eqref{Lorenz6},   with basis and coefficients in  Table \ref{tab:3}, we can also set up the partial differential equations for the learned mean residence time $u_{Learn}$,  and the learned escape probability $p_{Learn}$ together with the average escape probability $P_{Learn}$.

\begin{figure}[H]
\centering
\includegraphics[scale=0.15]{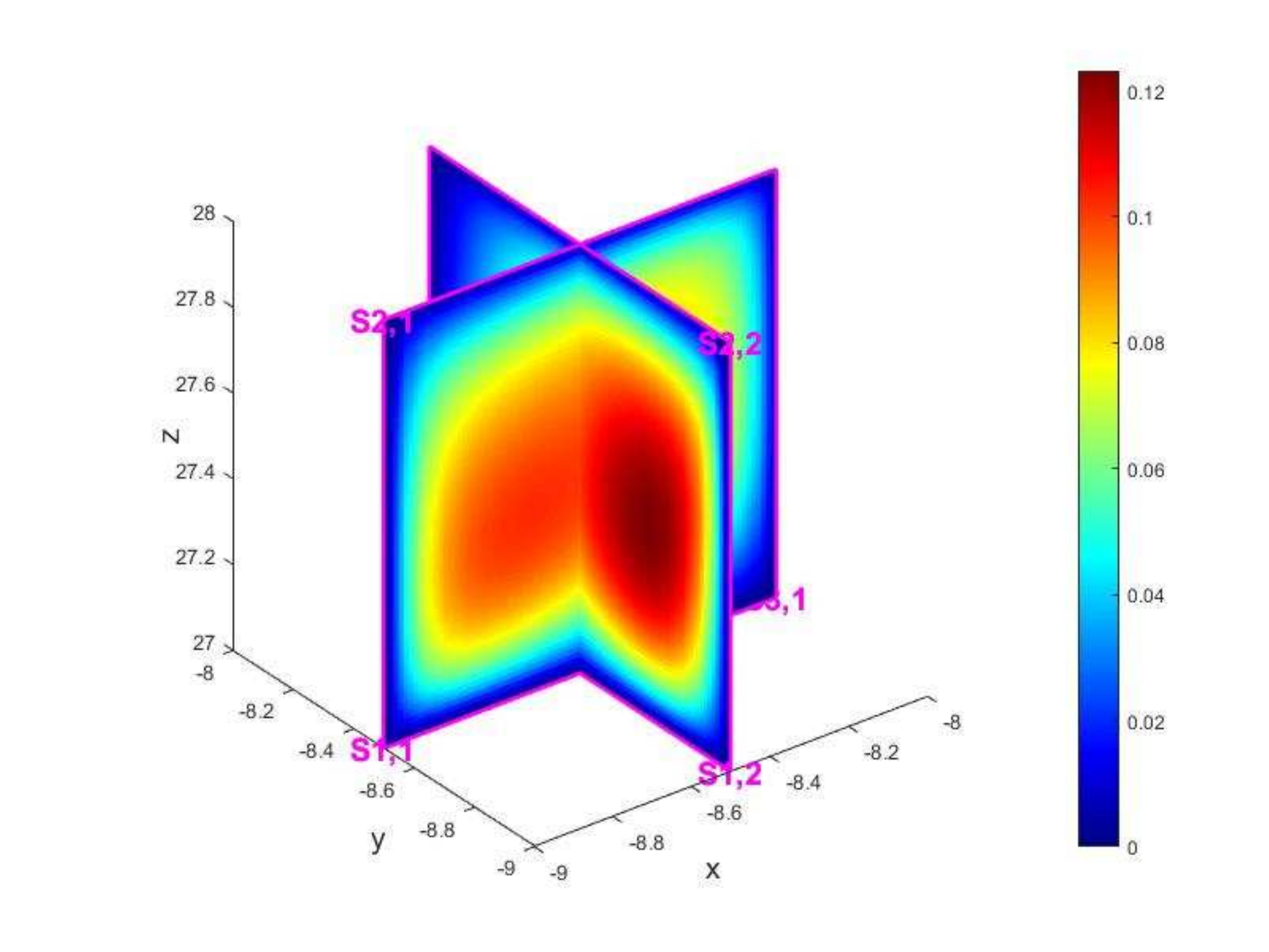}
\includegraphics[scale=0.15]{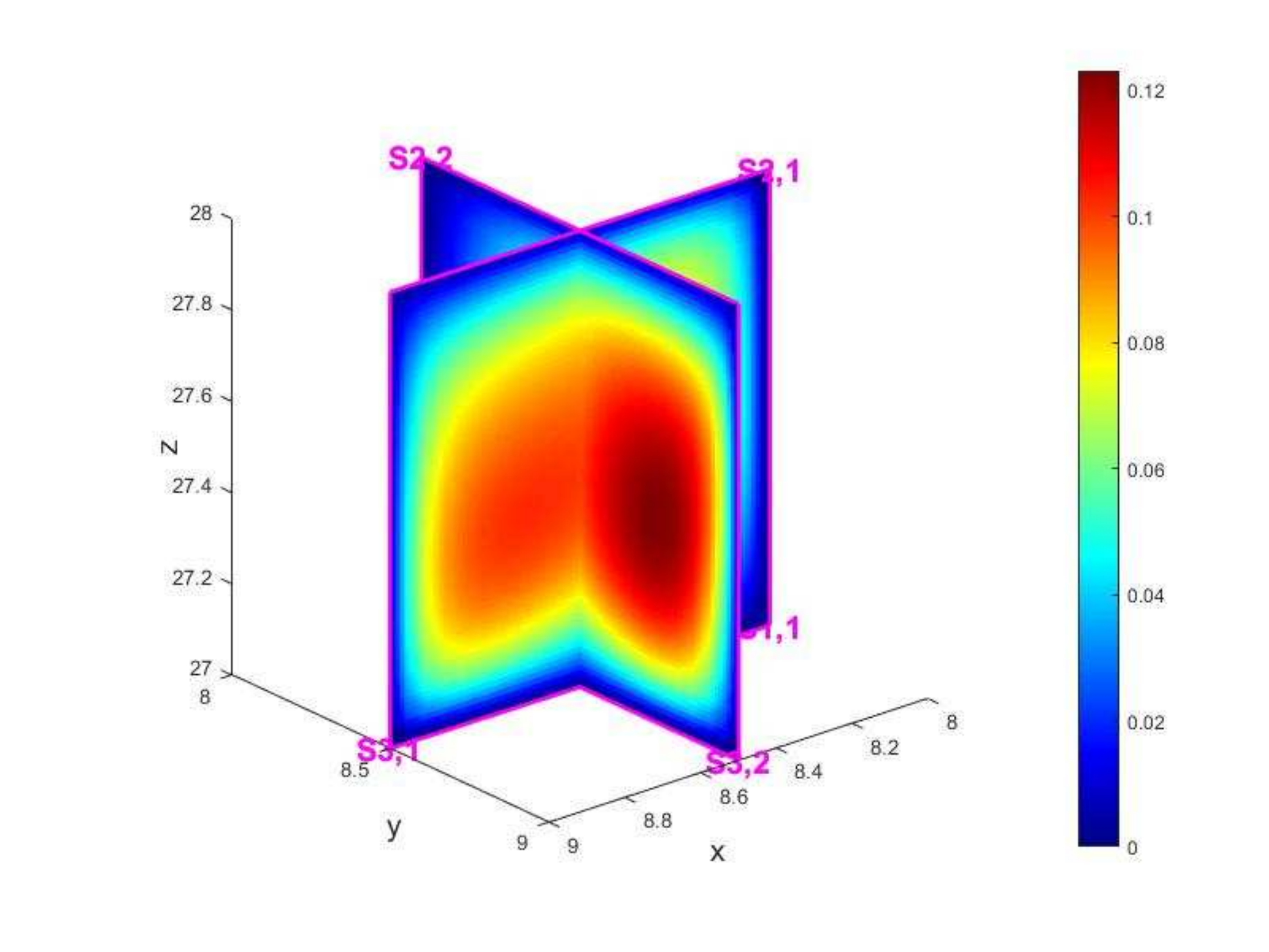}
\caption{Learned mean residence time of the learned Lorenz system \eqref{Lorenz4}-\eqref{Lorenz6}
 in $D_1$ (left) and $D_2$ (right): Two slices are shown for each subfigure.}
\label{figure 11}
\end{figure}

The learned mean residence time $u$ of the learned Lorenz
 system \eqref{Lorenz4}-\eqref{Lorenz6} is shown in FIG. \ref{figure 11}.
The maximal mean residence time in  the left region $D_1$ is $u_{Learn}(D_1)=0.1238$,
while the maximal mean residence time in the right region $D_2$ is
$u_{Learn}(D_2)=0.1236$. Then we compute  the error of maximal mean residence time between the learned  and the original  systems,
$\mbox{error}(u(D_1))=\max_{(x,y,z)\in D_1}\|u_{Learn} - u_{True}\|=1.0023\times 10^{-4}$,
and $\mbox{error}(u(D_2))=\max_{(x,y,z)\in D_{2}}\|u_{Learn} - u_{True}\|=1.2284\times 10^{-5}$.

\begin{figure*}[ht]
\centering
\includegraphics[scale=0.15]{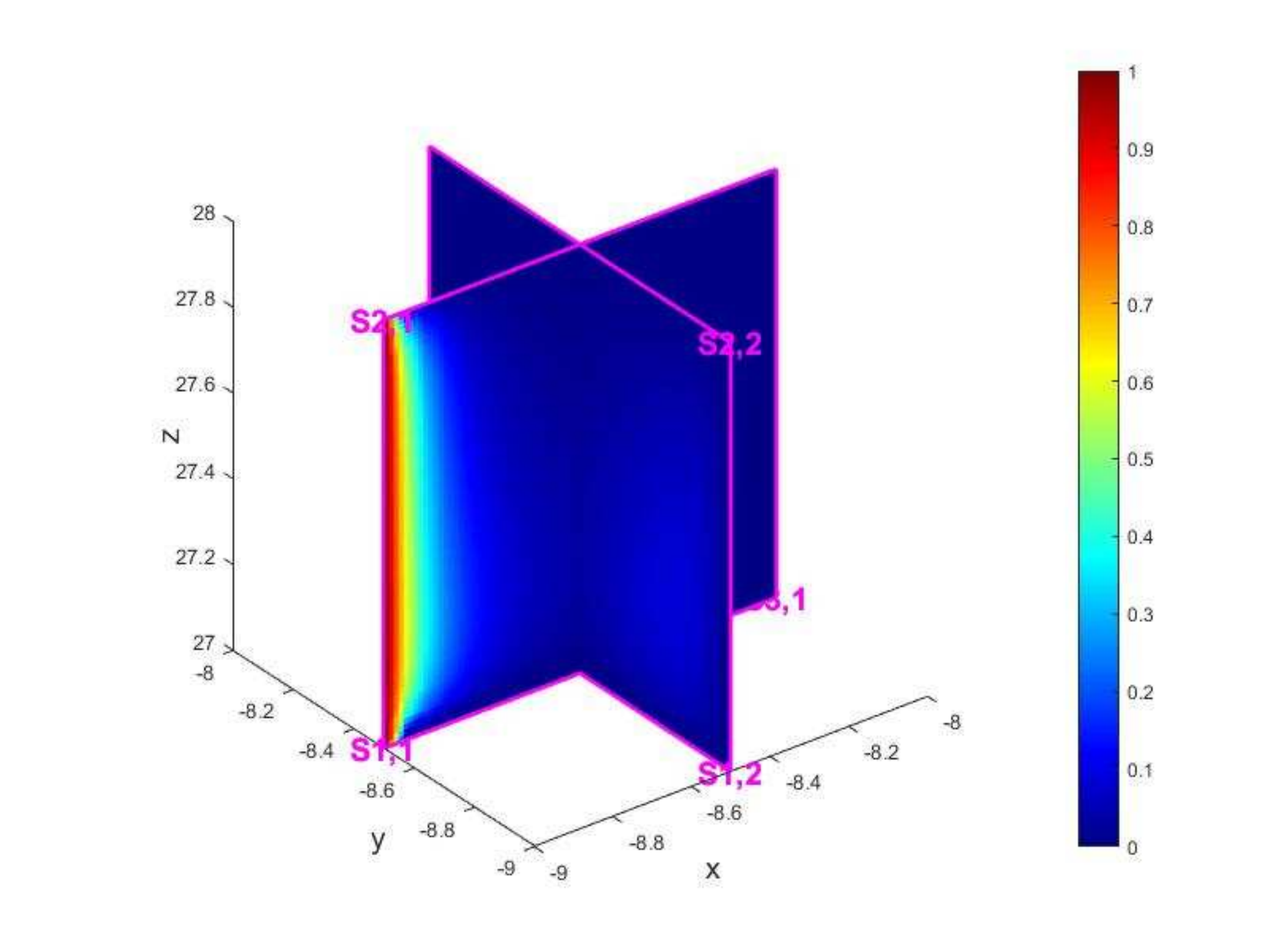}
\includegraphics[scale=0.15]{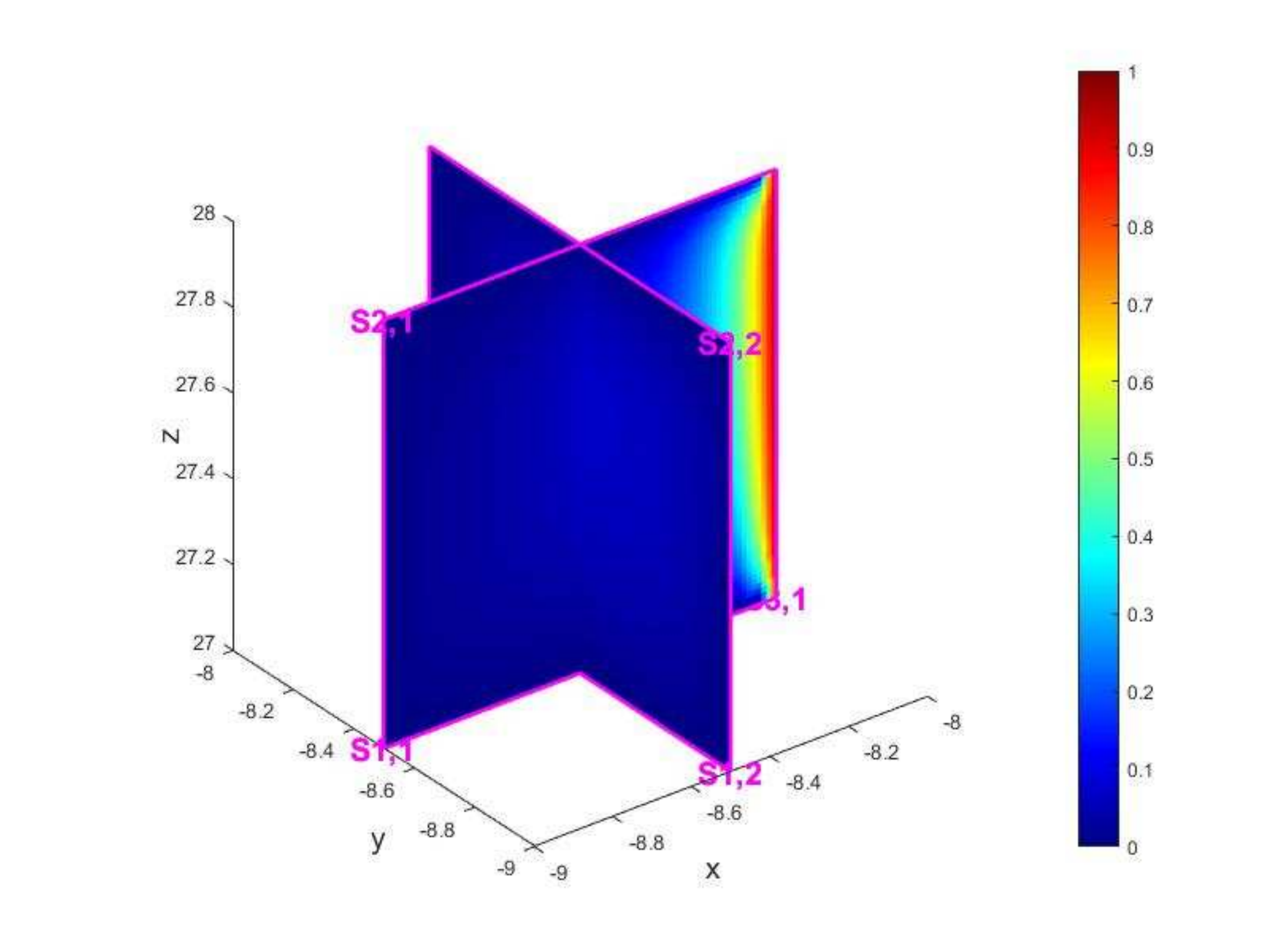}
\includegraphics[scale=0.15]{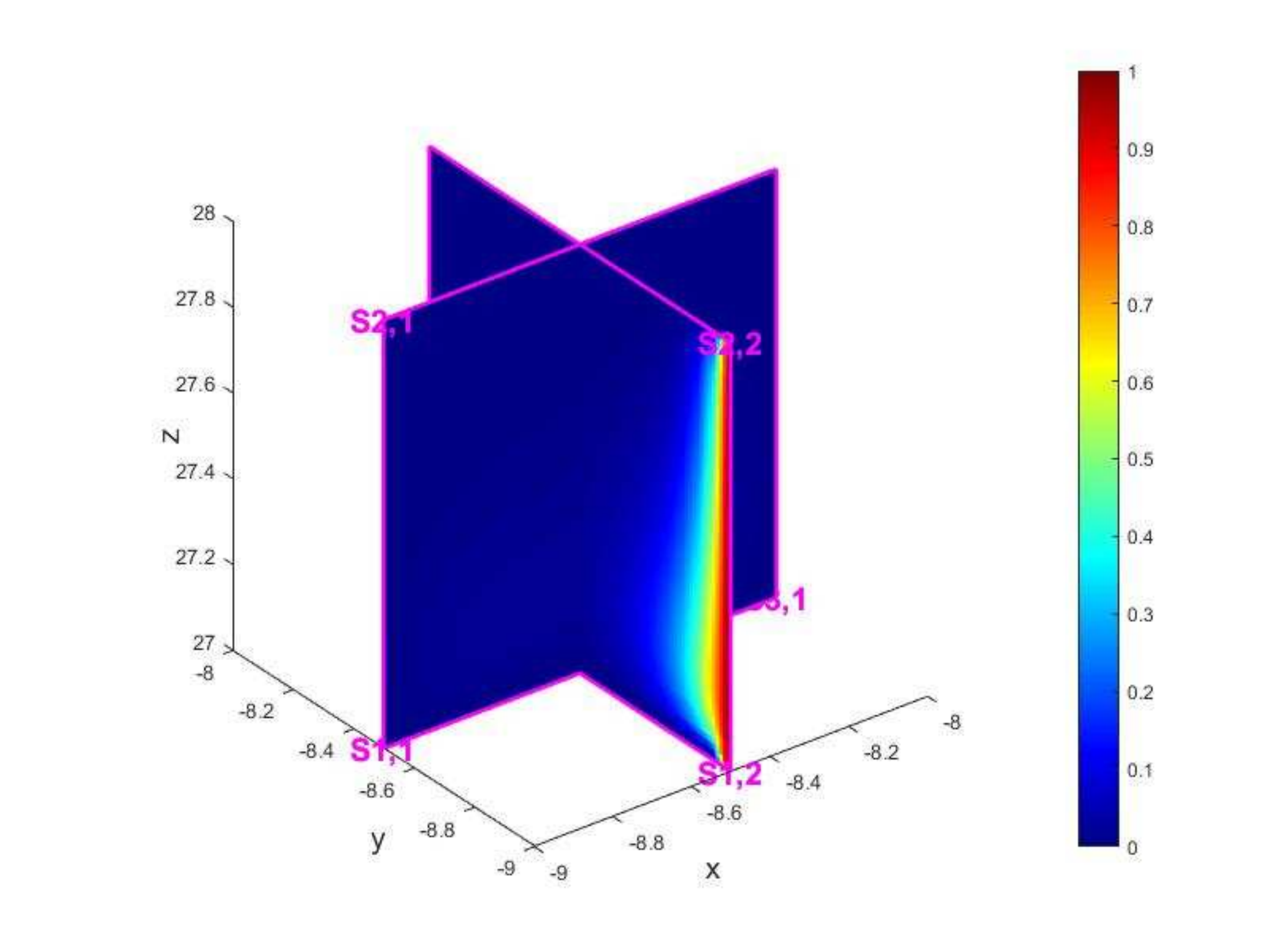}

\includegraphics[scale=0.15]{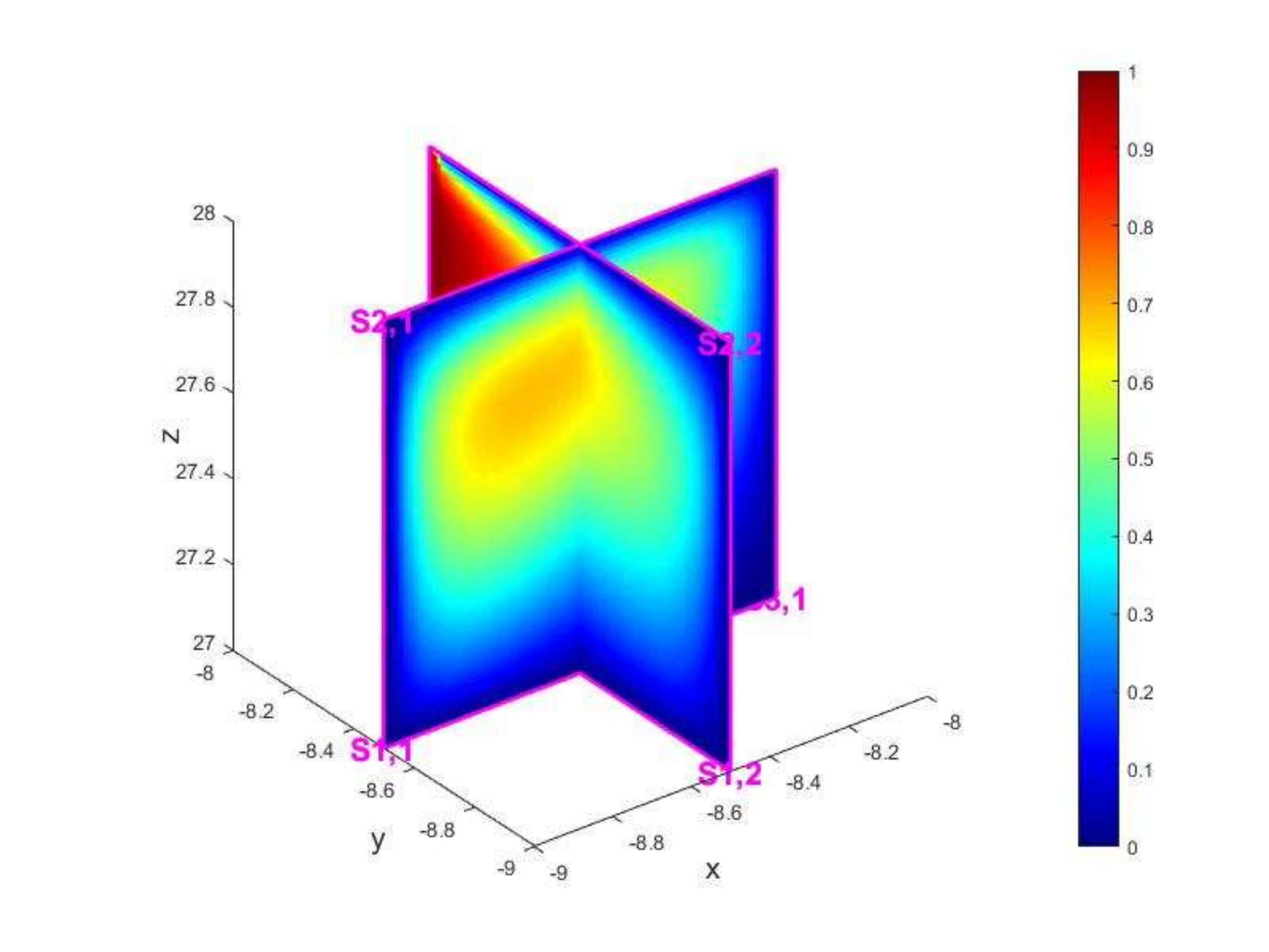}
\includegraphics[scale=0.15]{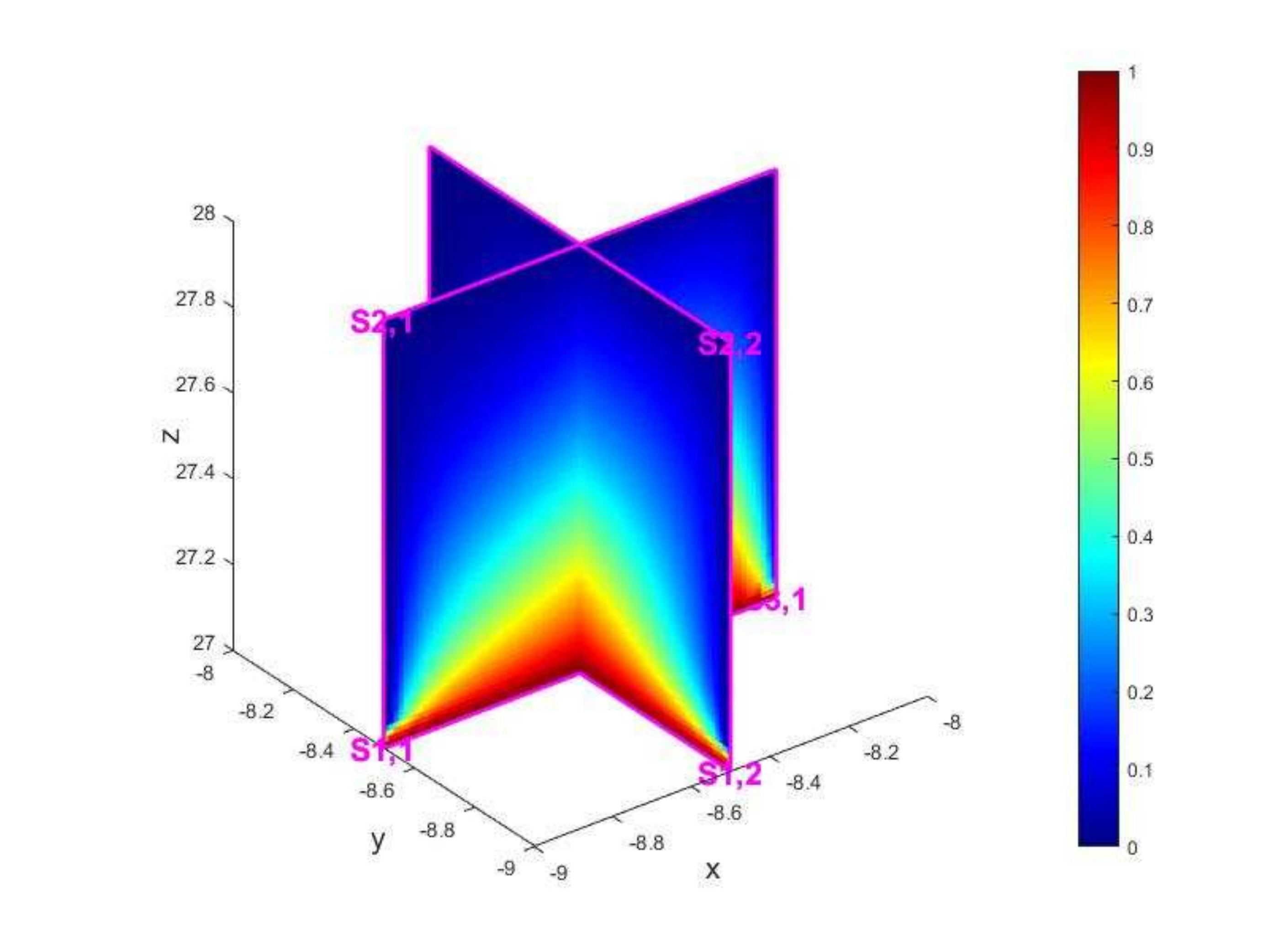}
\includegraphics[scale=0.15]{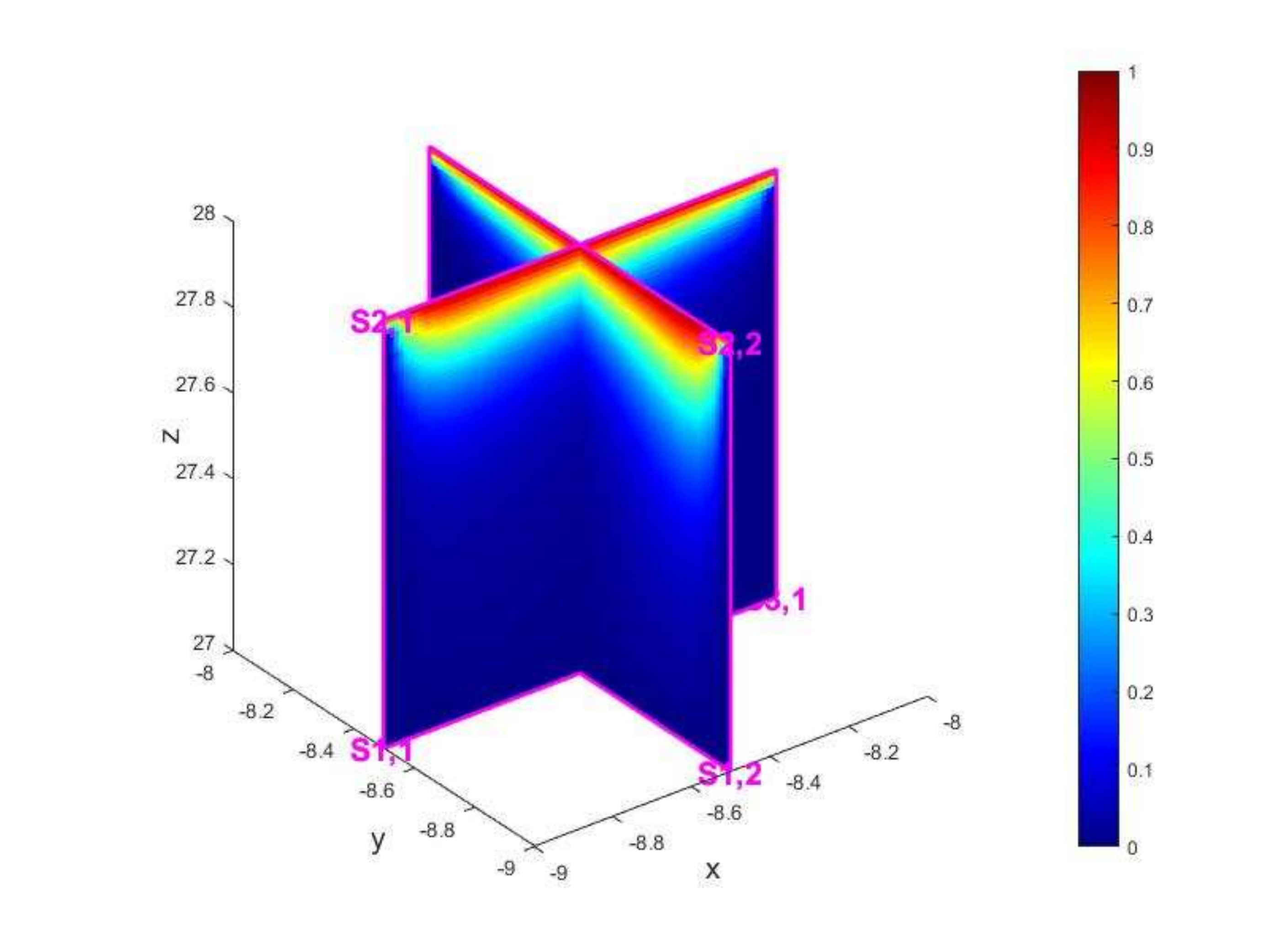}
\caption{Learned escape probability of the learned Lorenz system \eqref{Lorenz4}-\eqref{Lorenz6}  in
left   region $D_1$  exiting  from a subboundary $\Gamma$:   surface $x=-9$ (upper left),
$x=-8$ (upper middle), $y=-9$ (upper right), $y=-8$ (lower left),
$z=27$ (lower middle), and  $z=28$ (lower right).   Two slices are shown for each subfigure.}
\label{figure 12}
\end{figure*}

\begin{figure*}[ht]
\centering
\includegraphics[scale=0.15]{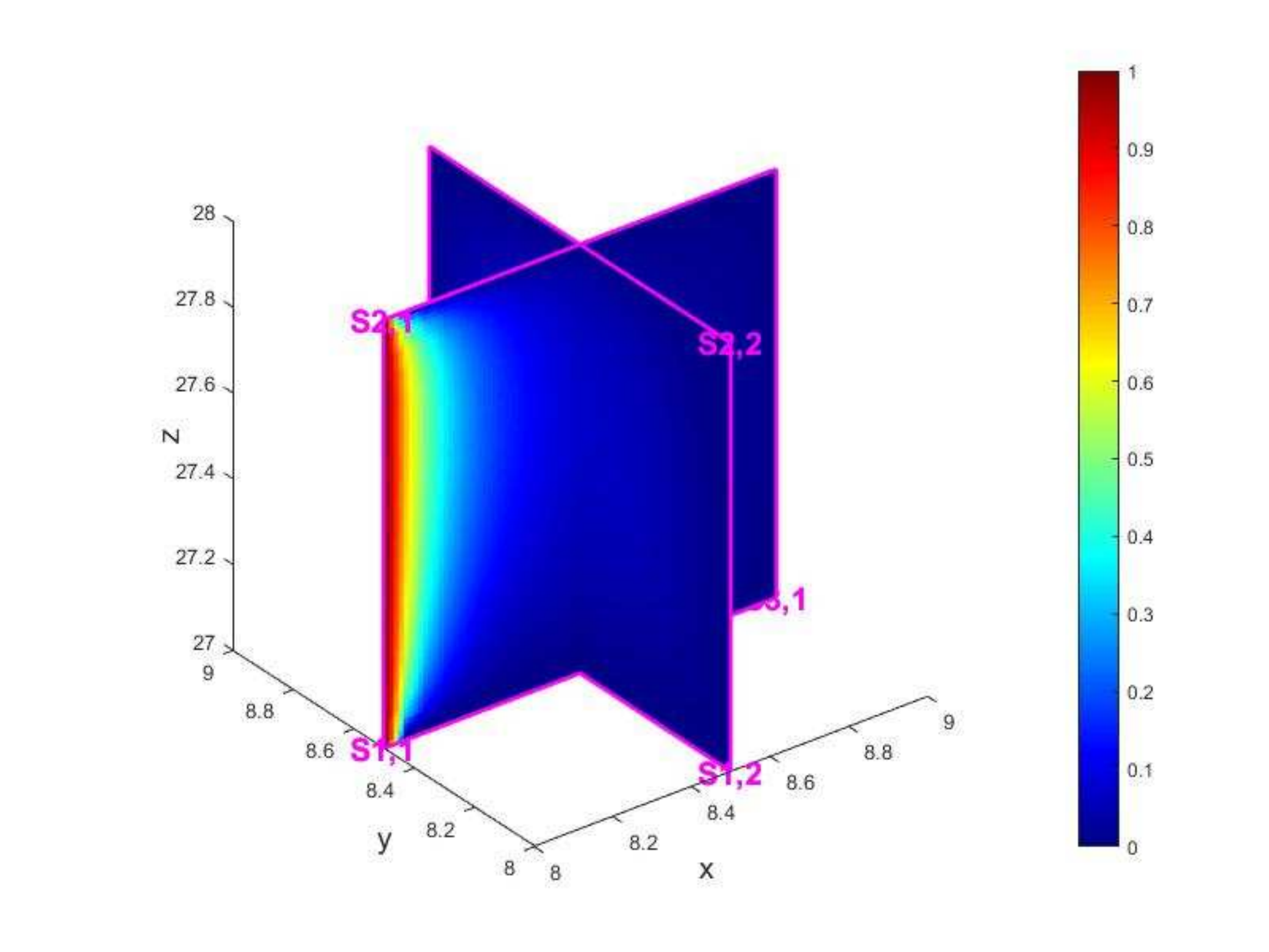}
\includegraphics[scale=0.15]{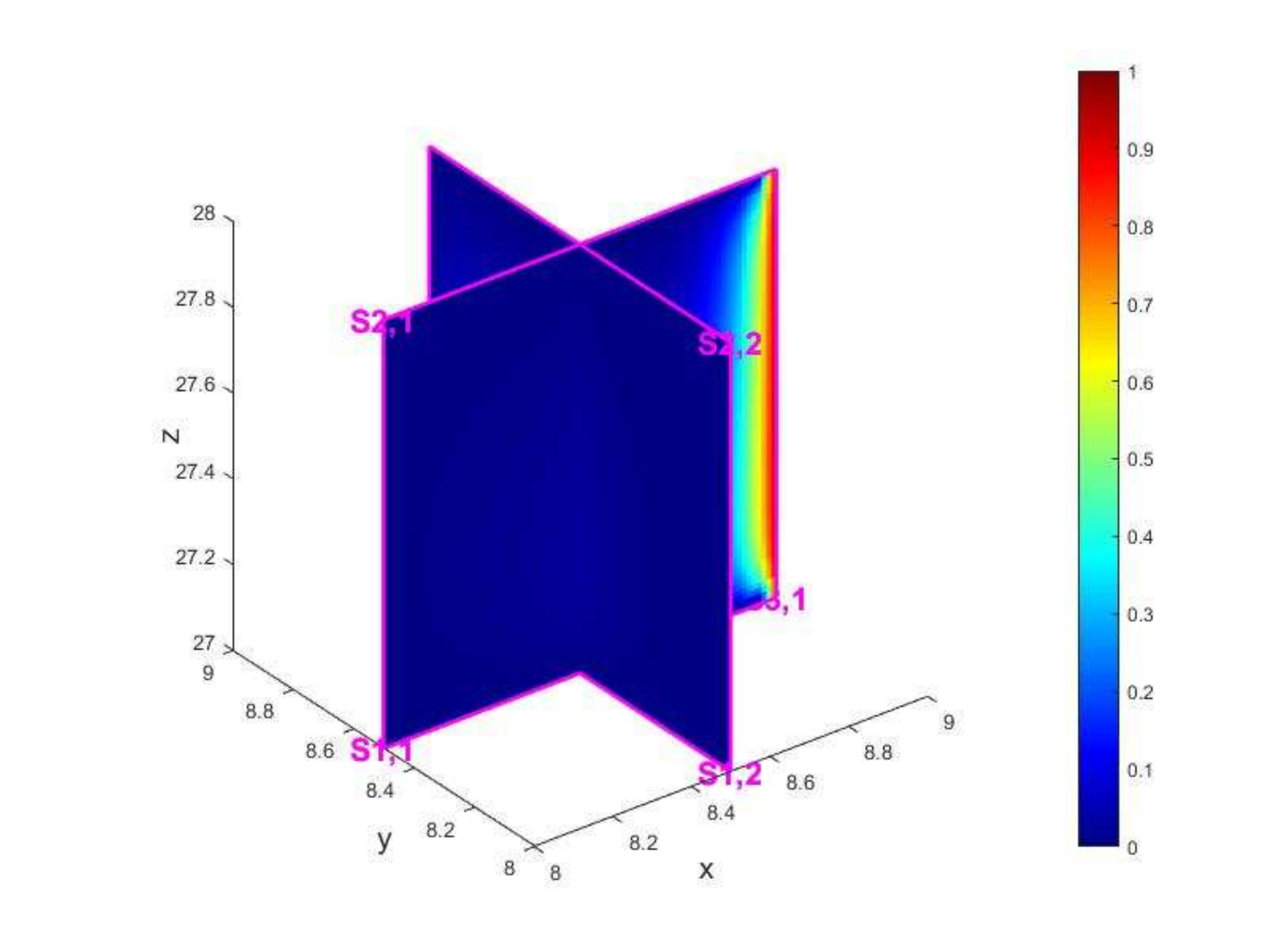}
\includegraphics[scale=0.15]{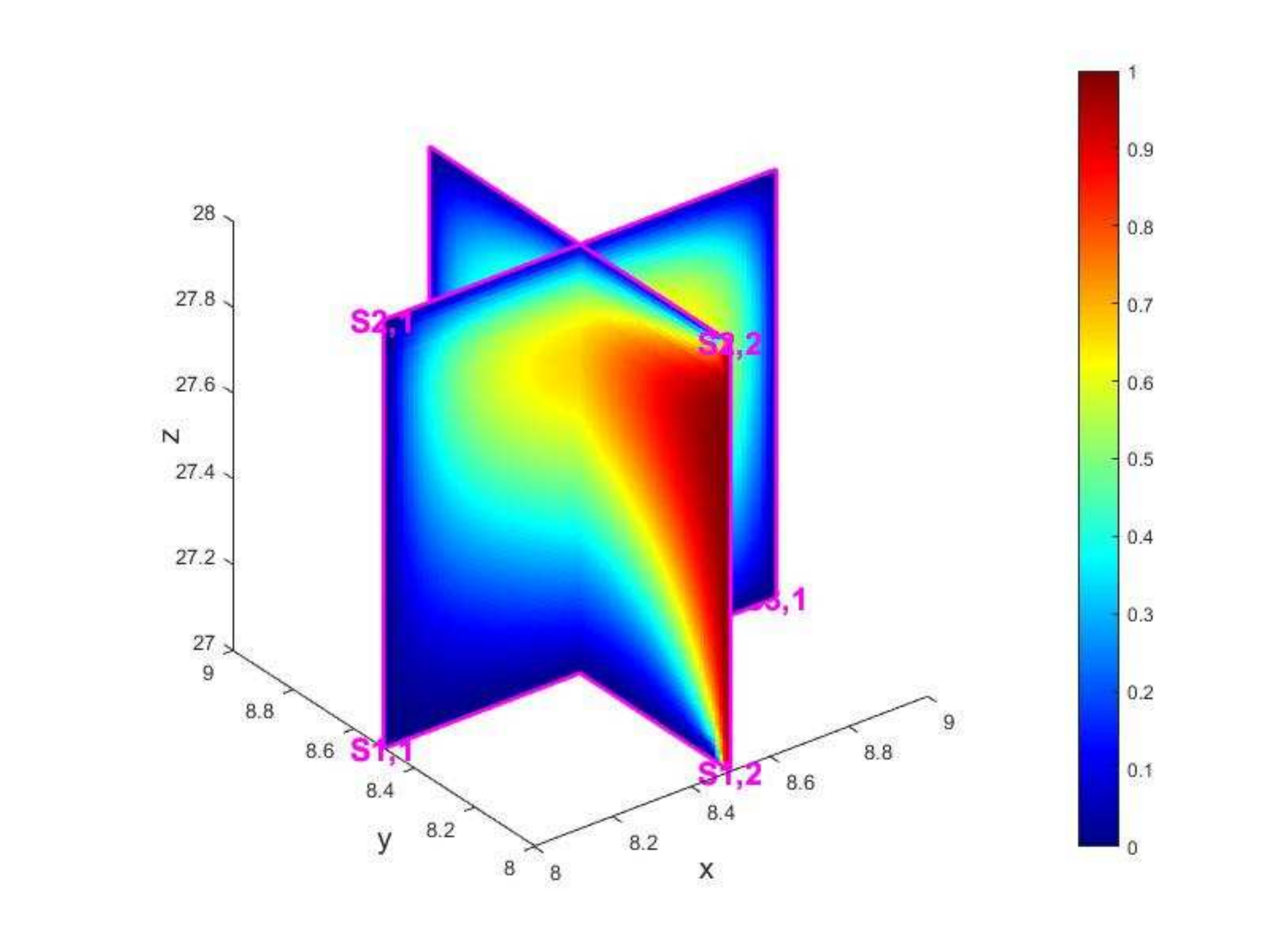}

\includegraphics[scale=0.15]{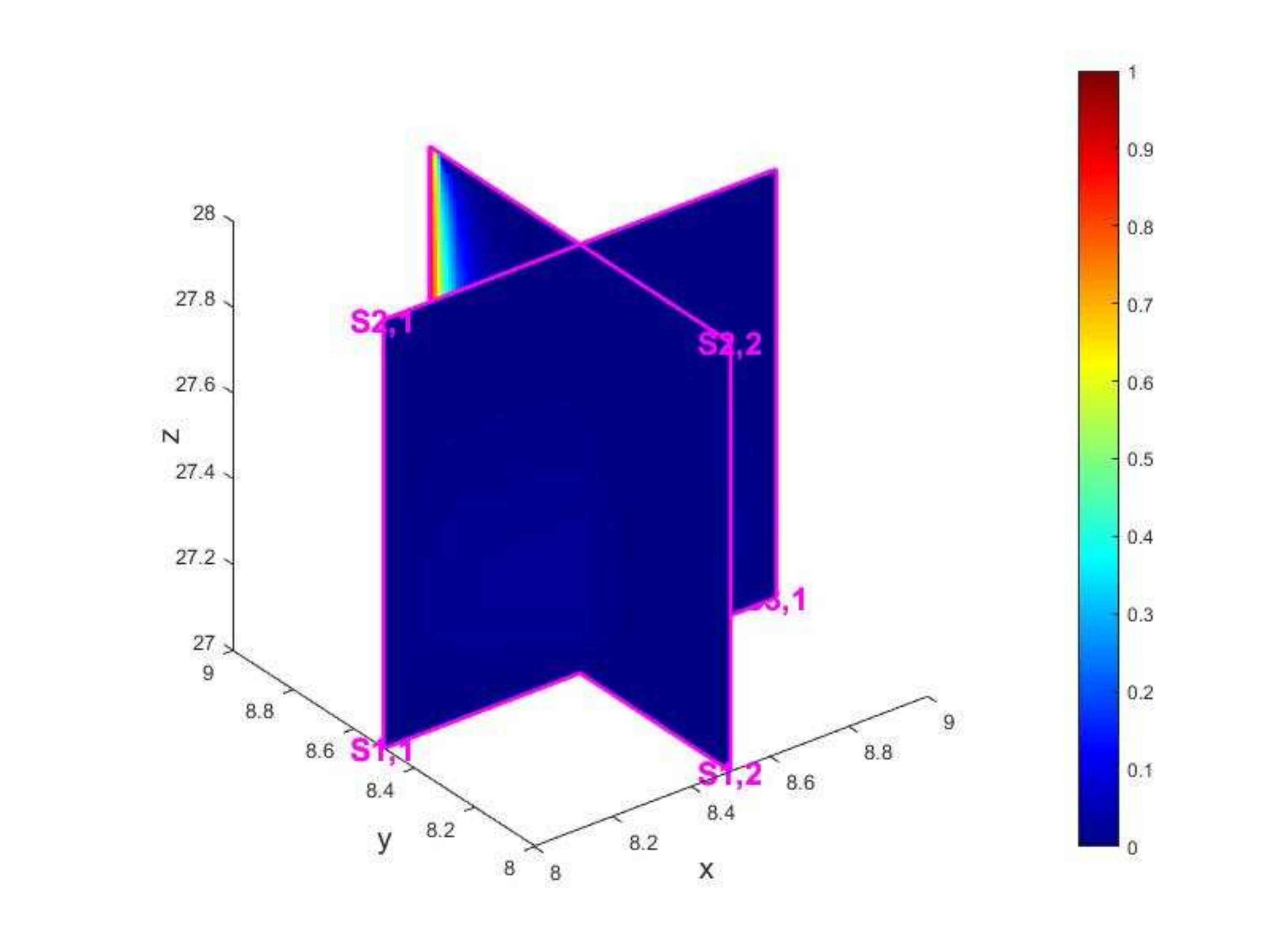}
\includegraphics[scale=0.15]{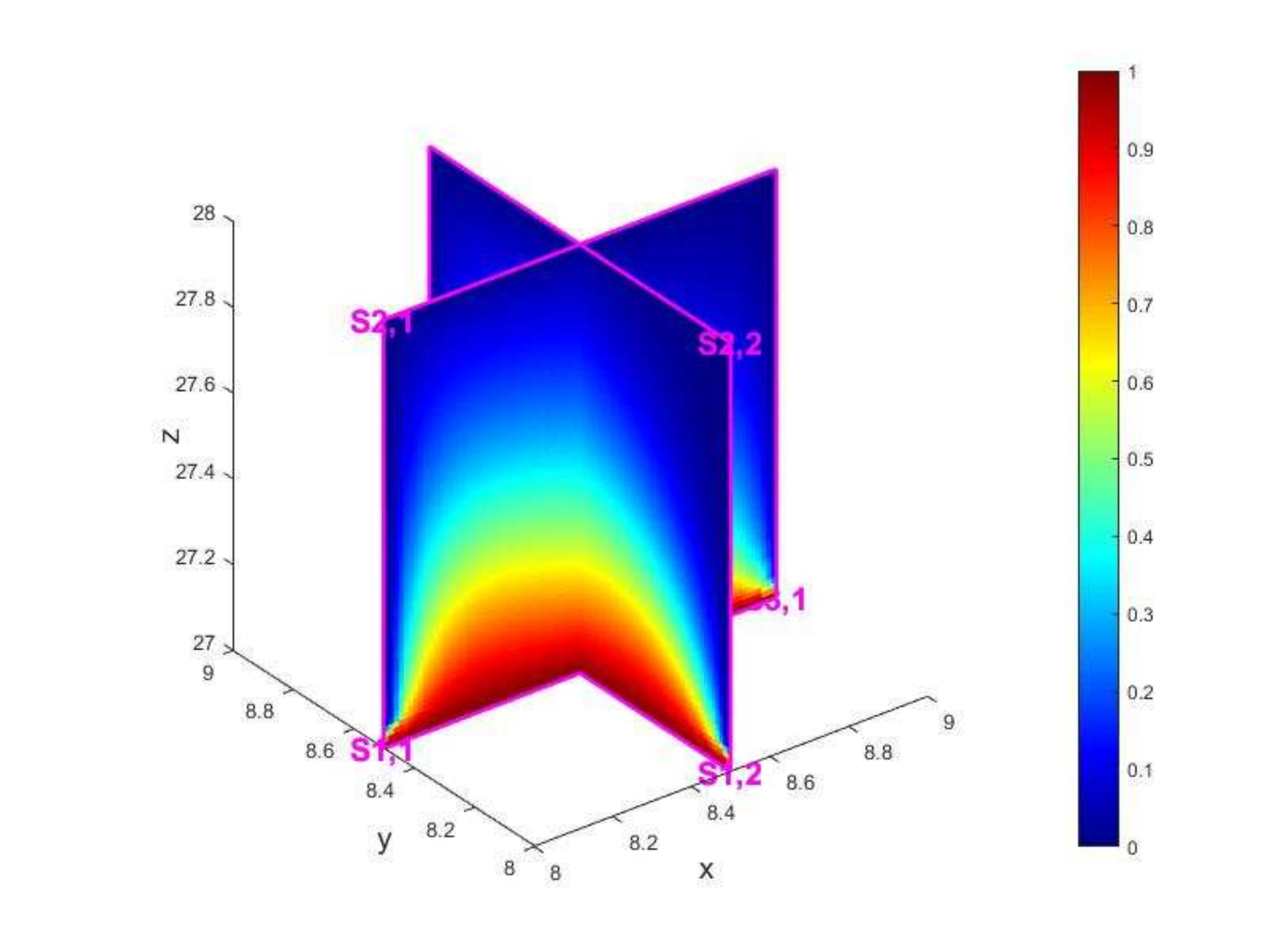}
\includegraphics[scale=0.15]{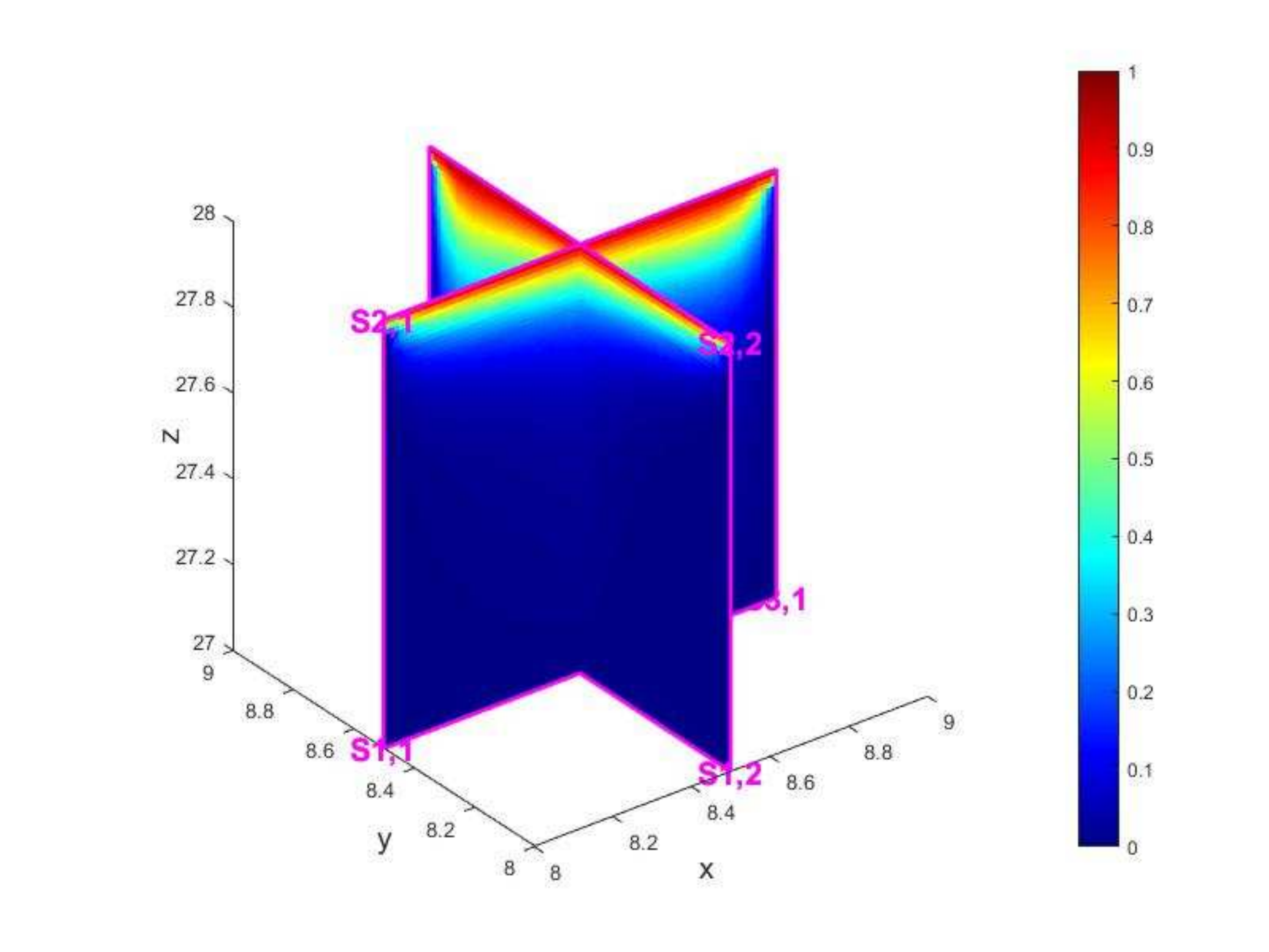}
\caption{Learned escape probability of the learned Lorenz system \eqref{Lorenz4}-\eqref{Lorenz6} in
right   region $D_2$  exiting from   a  subboundary $\Gamma$:  surface $x=8$ (upper left),
$x=9$ (upper middle), $y=8$ (upper right), $y=9$ (lower left),
$z=27$ (lower middle), and  $z=28$ (lower right).    Two slices are shown for each subfigure.          }
\label{figure 13}
\end{figure*}

We take $\Gamma$ to be each surface of the region.
The learned escape probability $p_{Learn}$ of particles exiting through each surface $\Gamma$ of left and right regions are  showed in
FIG. \ref{figure 12} - \ref{figure 13}, respectively, and
the average escape probability of left region is
$P_{Learn}(L_1) = 0.1065$ (from $x=-9$),
$P_{Learn}(L_2) = 0.1167$ (from $x=-8$),
$P_{Learn}(L_3)= 0.0955$ (from $y=-9$),
$P_{Learn}(L_4)= 0.3485$ (from $y=-8$),
$P_{Learn}(L_5)= 0.2343$ (from $z=27$),
$P_{Learn}(L_6)= 0.1244$ (from $z=28$).
Similarly, the average escape probability of right region is
$P_{Learn}(R_1)= 0.1161$ (from $x=8$),
$P_{Learn}(R_2)= 0.1061$ (from $x=9$),
$P_{Learn}(R_3) = 0.3469$ (from $y=8$),
$P_{Learn}(R_4)= 0.0947$ (from $y=9$),
$P_{Learn}(R_5)= 0.2367$ (from $z=27$),
$P_{Learn}(R_6) = 0.1264$ (from $z=28$).

Moreover, for example, we   compute the error     of  the  average escape probability between the learned  and the original  systems only for two cases,
 $\mbox{error}(P)(L4)=\|P_{Learn}(L4) - P_{True}(L4)\|=6.5838 \times 10^{-3}$,
 and $\mbox{error}(P)(R3)=\|P_{Learn}(R3) - P_{True}(R3) \|=6.431\times 10^{-3}$.
(L4)


\section{Conclusion}
\setcounter{equation}{0}
In summary, we have demonstrated a new data-driven method  to determine dynamical quantities,
mean residence time and escape probability,  based on sample path data of complex systems,  as long as these systems are  modeled by    stochastic   differential equations.

A novelty of this method is that it contains noise terms (in terms of Brownian motions) in the basis, in order to discover the governing stochastic differential equation. With the governing  stochastic differential equation as the model for the sample path data, we can then compute the mean residence time and escape probability, which are deterministic quantities that carry significant  dynamical  information.

This method combines machine learning tools (i.e., extracting governing stochastic differential equations from data) and stochastic dynamical systems techniques (i.e., understanding  random phenomena with deterministic quantities), to provide an efficient approach in examining stochastic dynamics from data.


To demonstrate that our method  is effective, we have illustrated  it on  several
  stochastic systems with Brownian motions.  It is expected to be  applicable to data sets from other stochastic dynamical systems.


\begin{acknowledgments}
We would like to thank Min Dai and Jian Ren for helpful comments.
The first author was supported by China Scholarship Council under No. 201808220008,
and the  second author was  supported by China Scholarship Council under No. 201808220009.
The second author  was also
supported by Jilin Province Industrial Technology Research and Development Special Project under No.2016C079.
The third author was partly supported by the National Science Foundation Grant 1620449.
\end{acknowledgments}

\bibliography{aipsamp}

\end{document}